\documentclass[reqno,11pt]{amsart} 
\usepackage{esint,mathrsfs,amssymb,amsmath}

\usepackage{bookmark,color}
\usepackage{comment}

\usepackage[left=1.2in, right=1.2in, bottom=1.5in]{geometry}
\newtheorem{theorem}{Theorem}[section]
\newtheorem{lemma}[theorem]{Lemma}
\newtheorem{proposition}[theorem]{Proposition}
\newtheorem{corollary}[theorem]{Corollary}
\theoremstyle{definition}
\newtheorem{definition}[theorem]{Definition}

\theoremstyle{remark}
\newtheorem{remark}[theorem]{Remark}

\numberwithin{equation}{section}

\newcommand{\norm}[1]{\lVert#1\rVert}
\newcommand{\Ag}[1]{\langle#1\rangle}

\newcommand{\R}{\mathbb{R}}
\newcommand{\Z}{\mathbb{Z}}

\newcommand{\T}{\mathbb{T}}
\newcommand{\e}{\varepsilon}
\newcommand{\txt}[1]{\text{\rm #1}}

\newcommand{\be}{\beta}

\newcommand{\vp}{\varphi}

\newcommand{\va}{\varepsilon}
\newcommand{\al}{\alpha}
\newcommand{\pa}{\partial}

\newcommand{\lm}{\lambda}
\newcommand{\Lm}{\Lambda}

\newcommand{\de}{\delta}

\newcommand{\na}{\nabla}

\newcommand{\X}{\mathcal{X}}

\newcommand{\vertiii}[1]{{\left\vert\kern-0.25ex\left\vert\kern-0.25ex\left\vert #1
		\right\vert\kern-0.25ex\right\vert\kern-0.25ex\right\vert}}

	\newcommand{\jz}[1]{\textcolor{red}{#1}}
	\newcommand{\fl}[1]{\lfloor #1 \rfloor}

\begin{document}

\title{Compactness and stable regularity in multiscale homogenization}



\author[W. Niu]{Weisheng Niu}
\address[W. Niu]{Center for Pure Mathematics, School of Mathematical Sciences, Anhui University, Hefei , China, 230601 P. R. China}
\email{niuwsh@ahu.edu.cn}

\author[J. Zhuge]{Jinping Zhuge}
\address[J. Zhuge]{Department of Mathematics, University of Chicago, Chicago, IL, 60637, USA}
\email{jpzhuge@uchicago.edu}



\begin{abstract}
	In this paper we develop some new techniques to study the multiscale elliptic equations in the form of $-\text{div} \big(A_\varepsilon \nabla u_{\varepsilon} \big) = 0$, where $A_\varepsilon(x) = A(x, x/\varepsilon_1,\cdots, x/\varepsilon_n)$ is an $n$-scale oscillating periodic coefficient matrix, and $(\e_i)_{1\le i\le n}$ are scale parameters. We show that the $C^\alpha$-H\"{o}lder continuity with any $\alpha\in (0,1)$ for the weak solutions is stable, namely, the constant in the estimate is uniform for arbitrary $(\varepsilon_1, \varepsilon_2, \cdots, \varepsilon_n) \in (0,1]^n$ and particularly is independent of the ratios between $\varepsilon_i$'s. The proof uses an upgraded method of compactness, involving a scale-reduction theorem by $H$-convergence. The Lipschitz estimate for arbitrary $(\varepsilon_i)_{1\le i\le n}$ still remains open. However, for special laminate structures, i.e., $A_\varepsilon(x) = A(x,x_1/\varepsilon_1, \cdots, x_d/\varepsilon_n)$, we show that the Lipschitz estimate is stable for arbitrary $(\varepsilon_1, \varepsilon_2, \cdots, \varepsilon_n) \in (0,1]^n$. This is proved by a technique of reperiodization.
\end{abstract}

\maketitle
\section{Introduction}
\subsection{Motivation}

This paper aims to study the uniform regularity of the multiscale elliptic equations (or systems) in the form of
\begin{equation}\label{eq.main}
	-\txt{div} \big(A_\e \nabla u_{\e} \big) = 0, \qquad \txt{in } \Omega \subset \R^d,
\end{equation}
where $A_\e(x) = A(x, x/\e_1,x/\e_2,\cdots, x/\e_n)$ is an $n$-scale oscillating coefficient matrix, and $(\e_i)_{1\le i\le n}$ are supposed to be small positive scale parameters. Throughout, we use the subscript $\e$ for $A_\e$ and $u_\e$ as we may think the sequence $(\e_i)_{1\le i\le n}$ relies implicitly on the variable parameter $\e>0$. Moreover, we assume that $A(x, y_1,y_2,\cdots,y_n)$ is 1-periodic with respect to every $y_i\in \R^d$. Without loss of generality (up to a reordering), we assume that $1\ge \e_1 \ge \e_2 \ge \cdots \ge \e_n>0$.

We investigate the uniform regularity of the solutions of \eqref{eq.main} that is independent of any $\e_i$.  For scalar equations, the classical De Giorgi-Nash estimate implies that the weak solution $u_\e$ of \eqref{eq.main} is $C^\alpha$-H\"{o}lder continuous for possibly some tiny $\alpha>0$ depending only on the ellipticity constant of $A_\e$ (not depending on the regularity or the structure of $A_\e$). For elliptic systems, such estimate does not exist, and the classical Schauder estimate will definitely depends on the regularity of $A_\e$ which is not uniform in any $\e_i$.
In the classical one-scale problem
(i.e., $n=1$ and $A_\e(x) = A(x/\e_1)$), quantitative homogenization and uniform regularity (up to Lipschitz estimate, if no corrector is involved) have been studied extensively in literatures in different settings; see monographs \cite{lions1978,Oleinik1988,ShenBook18}.

On the other hand, the periodic homogenization of \eqref{eq.main} with multiple scales ($n\ge 2$) has also been investigated, which turns out to be more difficult as the periodicity may be lost due to the coupling of different scales. In \cite{lions1978}, Bensoussan, Lions and Papanicolaou studied the case $\e_i = \e^i$ and showed that \eqref{eq.main} homogenizes to the same type of equation $-\txt{div} (\widehat{A} \nabla u_0) =0$ with some constant coefficient matrix $\widehat{A}$. The key insight in this case is that the fastest variable homogenizes first. Thus, we can find the homogenized operator on the finest microscale, next on the second finest, and continue repeatedly until all the microscales are homogenized. This process is called reiterated homogenization.  Later on, Allaire and Briane \cite{allaire1996} established the qualitative reiterated homogenization under a scale-separation condition, namely,
\begin{equation}\label{cond.separate}
	\lim_{\e\to 0} \e_1 = 0, \quad \lim_{\e \to 0} \frac{\e_{i+1}}{\e_i} = 0, \quad \txt{for all } i=1,2,\cdots, n-1.
\end{equation}
As far as the authors know, condition \eqref{cond.separate} was adopted in all the literatures in multiscale homogenization problems and is the minimal assumption to have a stable qualitative reiterated homogenization. To see this, if we consider $n=2$ and $\e_1 = \e, \e_2 = \lambda \e$ for some fixed constant $\lambda\in (0,1)$ (thus $\e_2/\e_1 = \lambda$), then the homogenized matrix will depend on $\lambda$ and the convergence rate is rather unstable, depending on the rationality or irrationality of $\lambda$ \cite{bbmm05}. However, if the condition \eqref{cond.separate} is assumed, then the homogenized matrix is fixed, and the optimal convergence rate could be  quantified \cite{nsxjfa2020}. See also \cite{past07} for related results in the case $n=2$ and $\e_2=(\e_1)^\alpha$.  Moreover, when the scales are well-separated, i.e., there exists $N\ge 1$ such that
\begin{equation}\label{cond.well-separate}
	\lim_{\e \to 0} \frac{1}{\e_i} \Big( \frac{\e_{i+1}}{\e_i} \Big)^{N} = 0, \quad \txt{for all } i = 1,2,\cdots, n-1,
\end{equation} 
 Niu, Shen and Xu \cite{nsxjfa2020} obtained the uniform Lipschitz regularity by a quantitative excess decay method. Condition \eqref{cond.well-separate} was first introduced by  Allaire and Briane in  \cite{allaire1996} to investigate qualitative reiterated homogenization of operators or perforated domains involving multiscales. This condition particularly includes the case $\e_i = \e^{\mu_i}$ with $0<\mu_1<\cdots< \mu_n < \infty$, but excludes the case $(\e_1, \e_2) = (\e, \e(|\ln \e| + 1)^{-1})$.  

A remaining fundamental question is that what happens if the scale-separation condition \eqref{cond.well-separate} or even \eqref{cond.separate} is not satisfied. Such situation arises more naturally in reality since all the oscillating scales (e.g., in composite materials) should be mesoscopic and do not go down to the scales of atoms or molecules. However, in this case the interaction of self-similar structures could be very difficult to handle quantitatively when two scales are close or comparable to each other, and the convergence rate could be extremely unstable (just as the example $\e_1 = \e$ and $\e_2 = \lambda \e$ when $\lambda$ varies between $0$ and $1$ \cite{bbmm05}). Or more generally, homogenization may fail for an arbitrary sequence of $(\e_i)_{1\le i\le n}$ that does not satisfy \eqref{cond.separate}. However, the failure of homogenization does not implies the non-existence of uniform regularity beyond De Giorgi-Nash estimate. In this paper, we will show surprisingly that nearly sharp uniform regularity holds stably for the equation \eqref{eq.main} for arbitrary $(\e_i)_{1\le i\le n}$ by an upgraded method of compactness. Moreover, in the special case that the coefficient matrix $A_\e$ is variable-separated (see \eqref{Var Separation}), we prove that the sharp uniform regularity estimates holds stably by a reperiodization argument.

\subsection{Assumptions and main results}
Throughout this paper, we denote by $d \ge 2$ the dimension and by $n \ge 1$ the number of scales under consideration. Let $A =   A(x,y_1,y_2,\cdots, y_n): \R^d\times \R^{d\times n} \mapsto \R^{d\times d}$ satisfy
\begin{itemize}
	\item Strong ellipticity condition: there exists $\Lambda \in (0,1]$ so that 
	\begin{equation}\label{ellipticity}
		\Lambda |\xi|^2 \le A(x,y_1,\cdots,y_n)\xi \cdot \xi \le \Lambda^{-1} |\xi|^2 
	\end{equation}for every $\xi \in \R^d$ and   for any $x \in \R^d, (y_1,\cdots,y_n)\in  \R^{d\times n}.$

	\item Periodicity: for any $(y_1,\cdots,y_n) \in \R^{d\times n}$ and $(z_1,\cdots,z_n)\in \Z^{d\times n}$,
	\begin{equation}\label{periodicity}
		A(x,y_1+z_1, \cdots, y_n+z_n) = A(x,y_1,\cdots,y_n).
	\end{equation}
	
	\item Smoothness: there exists $\gamma \in (0,1]$ and $L>0$ such that 
	\begin{equation}\label{smoothness}
		|A(x,y_1,\cdots,y_n) - A(x',y_1',\cdots,y_n')| \le L \big\{|x-x'|+  +\sum_{i=1}^n|y_i'-y_i|\big\}^\gamma 
	\end{equation} for any $x,x'\in \R^d$ and any $(y_1,\cdots,y_n), (y_1',\cdots,y_n')\in \R^{d\times n}$.
\end{itemize}

Consider the elliptic equation
\begin{equation}\label{eq.main.B1}
	-\txt{div} \big(A_\e \nabla u_{\e} \big) = 0 \qquad \txt{in } B_1 = B_1(0),
\end{equation}
where $A_\e(x) = A(x,x/\e_1,\cdots, x/\e_n)$ and $(\e_1,\e_2,\cdots, \e_n) \in (0,1]^n$.
The following is the first main result of this paper.
\begin{theorem}\label{thm.main1}
	Let $A$ satisfy \eqref{ellipticity} - \eqref{smoothness} and $u_\e\in H^1(B_1)$ be a weak solution of \eqref{eq.main.B1}. Then for every $\alpha \in (0,1)$, $u_\e$ is $C^\alpha$-H\"{o}lder continuous in $B_{1/2}$ uniformly in $(\e_i)_{1\le i\le n}$. Precisely, there exists $C_\alpha>0$ independent of $(\e_i)_{1\le i\le n}$, such that for any $x,y \in B_{1/2}$
	\begin{equation}\label{est.mainCa}
		|u_\e(x) - u_\e(y)| \le C_\alpha |x-y|^\alpha \norm{u_\e}_{L^2(B_1)}.
	\end{equation}
\end{theorem}

The originality of the above theorem is that $(\e_1,\e_2,\cdots, \e_n) \in (0,1]^n$ is arbitrary and the scale-separation condition is not needed. For this reason, we will say that the $C^\alpha$-H\"{o}lder regularity is stable in multiscale periodic homogenization. The notion of stability introduced here emphasizes that a (regularity) property is not only independent of the magnitudes of the scale parameters $(\e_i)_{1\le i\le d}$, but also independent of the ratios or any connections between them. Therefore, it is a stronger notion exclusively for multiscale homogenization, in place of the usual uniformity with respect to one scale parameter. 

We point out that \eqref{est.mainCa} combined with a blow-up argument implies that
\begin{equation}\label{est.almostLip}
    \norm{\nabla u_\e}_{L^\infty(B_{1/2})} \le C_\alpha \e_{n-1}^{\alpha-1}\norm{u_\e}_{L^2(B_1)}.
\end{equation}
When $\alpha = 1$, \eqref{est.mainCa} or \eqref{est.almostLip} is exactly the stable Lipschitz estimate, which is the optimal regularity one may expect for $u_\e$ if no correctors are involved. The above theorem shows that $u_\e$ has nearly optimal regularity in the sense that $\alpha$ can be arbitrarily close to $1$ (this could be compared to the one-scale situation in almost-periodic homogenization \cite{shenapde2015,armstrongcpam2016,shen-zhuge2016}). Unfortunately, the critical endpoint case $\alpha=1$ seems unapproachable through the proof in this paper; see Section 1.4 below.

Our next theorem deals with multiscale coefficients under a proper structural condition that admits the optimal stable Lipschitz estimate without condition (\ref{cond.separate}). Precisely, we say that $A_\e$ satisfies the variable-separation condition if
\begin{equation}\label{Var Separation}
	A_\e(x) = A(x,x_1/\e_1,x_2/\e_2,\cdots,x_d/\e_d),
\end{equation}
where $x = (x_1,x_2,\cdots, x_d) \in \R^d$.
Again, without loss of generality, we  assume $1\ge \e_1 \ge \e_2 \ge \cdots \ge \e_d>0$.

\begin{theorem}\label{thm.main2}
	 Assume $A = A(x,y): \Omega \times \R^d \to \R^{d\times d}$ satisfies the conditions (\ref{ellipticity})-(\ref{smoothness}).   Let $u_\e \in H^1(B_1)$ be a weak solution of $-\txt{div} (A_\e \nabla u_\e) = 0$ in $B_1$ with $A_\e(x)$ given by \eqref{Var Separation}. Then
	\begin{equation}\label{thm.main2-re}
		\norm{ \nabla u_\e}_{L^\infty(B_{1/2})} \le C\norm{u_\e}_{L^2(B_1)},
	\end{equation}
	where $C$ is independent of $(\e_i)_{1\le i\le n}$.
\end{theorem}

The structural condition of variable-separation arises naturally in multiscale laminates which play important roles in designing and manufacturing composite materials with certain properties; see \cite{tartar1985,milton1986,milton1999,chipot1986,briane2019,gloria2019} and references therein. In the special case of one-scale homogenization,  Dong, Li and Wang in \cite{dongdcds2018,dongjde2019} obtained the uniform interior Lipschitz estimate by the compactness method. In \cite{xuniucpde2020}, Xu and Niu established the optimal convergence rate and uniform regularity for more general one-scale laminate structure $A_\e = A(x, \rho(x)/\e)$, where $\rho:\Omega\to \R^n$ satisfies proper nondegenerate condition. Our result in Theorem \ref{thm.main2} provides a multiscale counterpart of the previous results, and also demonstrates that the Lipschitz estimate is stable, namely, it holds without separation of scales (\ref{cond.separate}). With the same idea, we generalize the Lipschitz estimate to spatially multiscale parabolic equations; see Subsection 4.2. Moreover, under the variable-separation condition we provide an effective approximation of the oscillating problem, and show that the effective approximate problem and the corresponding error estimate (convergence rate) are both stable; see Subsection 4.3. These stability properties are important in numerical computation for realistic problems.

We point out that the results in Theorem \ref{thm.main1} and Theorem \ref{thm.main2} can be generalized, without any real difficulties, to boundary estimates, elliptic systems, or having a source term on the RHS, etc. 

\subsection{Main ingredients of the proofs}
We first discuss the proof of Theorem \ref{thm.main1}. Since \eqref{cond.well-separate} is not assumed in Theorem \ref{thm.main1}, the quantitative method developed recently \cite{armstrongan2016,shenapde2017} does not apply. Instead, we develop an upgraded method of compactness, originating from Avellenda and Lin \cite{al87,al89}. However, our situation is logically more complicated and tortuous than the one-scale problems, and essentially needs some new ideas to break the barrier, which will be explained below briefly. First, the compactness method relies only on the qualitative homogenization, usually called $H$-convergence, introduced by Tartar (see Definition \ref{H-convergence}). The main step in our argument is a scale-reduction process in Theorem \ref{thm.ScaleReduction}. Precisely, we will show that for any given sequence of $(\e_1^{(k)},\e_2^{(k)},\cdots, \e_n^{(k)})_{k\ge 1}$ with $\lim_{k\to 0} \e_n^{(k)} = 0$, there is a subsequence and some integer $m$ with $0\le m< n$ so that
\begin{equation}\label{eq.H-conv}
	A(x, x/\e^{(k_j)}_{1},\cdots,x/\e^{(k_j)}_{n}) \quad H\txt{-converges to}\quad B(x, x/\delta_1,\cdots, x/\delta_m),
\end{equation}
where $B$ is a coefficient matrix satisfying \eqref{ellipticity} - \eqref{smoothness} and $(\delta_1, \delta_2,\cdots, \delta_m) \in (0,\infty)^m$ ($B = B(x)$ if $m=0$). In other words, an $n$-scale periodic coefficient matrix can be reduced to an $m$-scale periodic coefficient matrix through $H$-convergence with $m$ strictly less than $n$. Hence, through finite times of scale-reduction, we may arrive at a one-scale problem whose regularity is well-understood. This allows us to prove the uniform $C^\alpha$ regularity by an induction on the number of scales $n$.

The proof of \eqref{eq.H-conv} is a delicate rearrangement of the scales and variables combined with reiterated homogenization of quasi-periodic coefficient matrices (see Theorem \ref{thm.quasiper.HC}). To demonstrate the key idea, we consider a simple case $A_\e(x) = A(x/\e_1, x/\e_2)$ and assume $\e_2 \le \e_1 $ and $\lim_{\e \to 0} \e_1 = \lim_{\e \to 0} \e_2 = 0$. If $\lim_{\e \to 0} \e_2/\e_1 = 0$ (up to a subsequence), then $A_\e$ $H$-converges to a constant matrix by reiterated homogenization \cite{allaire1996}. The difficulty lies in the case $\lim_{\e \to 0} \e_2/\e_1 = \lambda > 0$ (up to a subsequence). The key insight is that, in this case, we may rewrite
\begin{equation}\label{eq.Ae.n=2}
	A_\e(x) = A\Big(\frac{x}{\e_1}, \Big( \frac{\e_1}{\e_2} - \frac{1}{\lambda}\Big) \frac{x}{\e_1} + \frac{x}{\lambda \e_1}  \Big),
\end{equation}
and define
\begin{equation}\label{eq.e'A'.n=2}
	\e_2' = ( \frac{\e_1}{\e_2} - \frac{1}{\lambda})^{-1} \e_1 \quad \txt{and} \quad A'(y_1, y_2) = A(y_1, y_2 + \lambda^{-1}y_1).
\end{equation}
Then the rearranged matrix $A'(y_1, y_2)$ is quasi-periodic in $y_1$ and periodic in $y_2$. Moreover, observe that $A_\e(x) = A'(x/\e_1, x/\e_2')$ and $\lim_{\e \to 0} \e_1/\e_2' = 0$. This means that $A'(x/\e_1, x/\e_2')$ is quasi-periodic and satisfies the scale-separation condition \eqref{cond.separate}. Now, depending on the limit $\lim_{\e \to 0} \e_2'$, we have three subcases. If $\lim_{\e \to 0} \e_2' = 0$, then Theorem \ref{thm.quasiper.HC} shows that $A_\e$ $H$-converges to a constant matrix; if $\lim_{\e \to 0} \e_2' = \delta\in (0,\infty)$, then $A_\e$ $H$-converges a periodic matrix in a form of $B(x/\delta)$ due to Lemma \ref{lem.AjBj}; finally if $\lim_{\e \to 0} \e_2' = \infty$, then $A_\e$ also $H$-converges to a constant matrix due to Lemma \ref{lem.AjBj}. In any case above, the number of scales has been reduced from two to one or zero. The above argument is a simple case with only two scales. For the most general case with arbitrary sequence of $(\e_1^{(k)},\e_2^{(k)},\cdots, \e_n^{(k)})_{k\ge 1}$, the argument is much more subtle as we need to make sure that all the scales are separated by a careful rearrangement in spirit of \eqref{eq.Ae.n=2} and \eqref{eq.e'A'.n=2}. We mention that Theorem \ref{thm.quasiper.HC}, regarding the reiterated homogenization for quasi-periodic coefficients, was claimed in \cite{guenuean2019} without a proof. We provide a proof in Section 5 by generalizing the method of periodic multiscale convergence \cite{allaire1996}, which is of independent interest.

Next, we discuss briefly the proof of Theorem \ref{thm.main2}. We actually provide two different proofs of Theorem \ref{thm.main2}. The first one is called a method of reperiodization which seems new and rather simple. The idea is to take advantage of the laminate structure of the coefficients and perform a suitable change of variables to reperiodize the coefficient matrix (finitely many times) so that the existing result may apply directly and repeatedly. A simple example of reperiodization may be seen in Subsection \ref{sec.open equation}.
 This method can be adapted easily to parabolic equations with non-self-similar scales (for time-dependent coefficients) and gives simple proofs to the regularity estimates contained in \cite{GS20,GN21}; see Theorem \ref{thm.parabolic}.

The second approach (also motivated by the recent work \cite{GS20,GN21}) is based on the excess decay method developed recently in \cite{armstrongan2016,shenapde2017}. One of the key step is to find proper smooth functions to approximate the solution $u_\va$ of the original oscillating problem. And generally, the homogenized solution would be a good choice, see \cite{armstrongan2016,shenapde2017}. Unfortunately, when the scale-separation condition \eqref{cond.separate} is not satisfied, there is no stable homogenized solution with desirable error estimate. Our idea here is to approximate $u_\va$ by a family of effective solutions (solutions of the effective approximate problems) depending on $(\e_i)_{1\le i\le d}$. The key point is then to show that the convergence rate to the effective solutions is stable with respect to $(\e_i)_{1\le i\le d}$; see Theorem \ref{error-th}. This requires us to obtain some uniform estimates related to certain degenerate cell problems, which can be handled by the idea of \cite{GS20,GN21} or the idea of reperiodization mentioned above.

\subsection{An open question}\label{sec.open equation}
As we have mentioned, Theorem \ref{thm.main1} does not include the Lipschitz estimate, namely, the endpoint case $\alpha = 1$ in \eqref{est.mainCa}. On the other hand, we know that the Lipschitz estimate is true at least for two special cases: (1) the scales are well-separated \cite{nsxjfa2020}, or (2) under the variable-separation condition in Theorem \ref{thm.main2}. It is critical to point out that besides the above two cases, we can still find other special cases that admit Lipschitz estimate. For example (this is a simple example of reperiodization), consider $n = 2$ and $\e_1 = m\e, \e_2 = \e$ with $m\ge 1$ being a fixed integer. In this case, $A_\e(x) = A(x/(m\e), x/\e)$ and $\e_2/\e_1 = 1/m$. The key observation is that if we set $\tilde{A}(y) = A(y,my)$, then $\tilde{A}(y)$ is 1-periodic in $y$ and $\tilde{A}(x/(m\e)) = A_\e(x)$. Hence, by the well-known large-scale Lipschitz estimate for the one-scale problem, we have
\begin{equation}\label{est.Am.Lip}
	\sup_{m\e < r<1} \bigg( \fint_{B_r} |\nabla u_\e|^2 \bigg)^{1/2} \le C\bigg( \fint_{B_1} |\nabla u_\e|^2 \bigg)^{1/2}.
\end{equation}
Next, for $0<r<m\e$, we apply a blow-up argument to reach a coefficient matrix in a form of $A(x,mx)$. This is a standard problem with only one periodically oscillating scale, whose Lipschitz estimate is also known \cite{nsxjfa2020}. This, combined with \eqref{est.Am.Lip}, gives the pointwise Lipschitz estimate uniform in both $\e$ and $m$. Note that the above example is extremely sensitive on $m$ and the argument relies essentially on the fact that $m$ is an integer which allows reperiodization. If $m$ is not an integer, we may also expect a good Lipschitz estimate if $m$ is a suitable rational number or an irrational number satisfying certain Diophantine condition \cite{armstrongcpam2016,agk16,shen-zhuge2016}. Unfortunately, all these special cases do not conclude the stable Lipschitz estimate in general setting (even for $n=2$) of Theorem \ref{thm.main1}. We now restate this unsolved problem as an open question.

\textbf{Open question:} Under the assumptions of Theorem \ref{thm.main1}, does \eqref{est.mainCa} hold for $\alpha = 1$? Or in other words, is Lipschitz estimate stable in multiscale homogenization?

As a weaker question, one may also ask for the $W^{1,p}$ estimate for a large $p>2$, which is also unclear.

\subsection{Organization of the paper}
In Section 2, we study the compactness under $H$-convergence and prove a scale-reduction theorem. This theorem then is used in Section 3 to prove Theorem \ref{thm.main1}. In Section 4, we consider the case of variable-separation. Particularly, we prove Theorem \ref{thm.main2} in Subsection 4.1 and then apply our method to parabolic equations in Subsection 4.2. The {existence and} stability of the effective approximation is discussed in Section 4.3. In Section 5, we prove the reiterated homogenization theorem for quasi-periodic coefficients, i.e., Theorem \ref{thm.quasiper.HC}.

\subsection*{Acknowledgements}
W. Niu is supported by NNSF of China (11971031) and  NSF of Anhui Province (2108085Y01)

\section{$H$-convergence and scale reduction}

In this section, we study the compactness of $H$-convergence and prove a theorem of scale reduction. Let us first recall the definition of $H$-convergence.

Let $\Omega\subset \R^d$ be a bounded domain. For $0<\alpha\le \beta < \infty$, let $\mathcal{M}(\alpha,\beta;\Omega)$ denote the set of $A\in L^\infty(\Omega;\R^{d\times d})$ satisfying
\begin{equation*}
	A^{-1}(x)\xi\cdot \xi \ge \beta^{-1}|\xi|^2, \quad A(x)\xi\cdot \xi \ge \alpha |\xi|^2, \quad \text{a.e. } x\in \Omega.
\end{equation*}

\begin{definition}[\cite{Tbook09}]\label{H-convergence}
	A sequence $(A_k)_{k\ge 1} \subset \mathcal{M}(\alpha,\beta;\Omega)$ $H$-converges to $\bar{A}\in \mathcal{M}(\alpha',\beta';\Omega)$ for some $0<\alpha'\le \beta'<\infty$, if for all $f\in H^{-1}(\Omega)$, the sequence of solutions $(u_k)_{k\ge 1} \subset H^1_0(\Omega)$ of $-\txt{div} (A_k \nabla u_k) = f$ converges to $\bar{u}$ weakly in $H^1_0(\Omega)$; and the sequence $A_k \nabla u_k$ converges to $\bar{A} \nabla \bar{u}$ weakly in $L^2(\Omega;\R^d)$, where $\bar{u}$ is the weak solution of $-\txt{div} (\bar{A} \nabla \bar{u}) = f$ in $\Omega$.
\end{definition}

The reason of considering the class $\mathcal{M}(\alpha,\beta;\Omega)$ in $H$-convergence is that it is compact under $H$-convergence \cite[Theorem 6.5]{Tbook09}. In other words, for any sequence $(A_k)_{k\ge 1} \subset \mathcal{M}(\alpha,\beta;\Omega)$, there exists a subsequence that $H$-converges to some $\bar{A}\in \mathcal{M}(\alpha,\beta;\Omega)$. Also note that a matrix satisfying \eqref{ellipticity} belongs to $\mathcal{M}(\alpha,\beta;\Omega)$ for some $\alpha,\beta>0$ depending only on $\Lambda$.

In the following, we assume that $(\e_1, \e_2,\cdots, \e_n)$ satisfies
\begin{equation}\label{cond.order}
	\infty> \e_1\ge \e_2 \ge \cdots \ge \e_n > 0.
\end{equation}
Even though previously we assume $\e_i \le 1$, the case $\e_i > 1$ will take place under rescaling which is crucial in our argument.
Let $\e^{(k)} = (\e^{(k)}_{1}, \e^{(k)}_{2},\cdots, \e^{(k)}_{n})$, $k\in 1,2,\cdots$, be a sequence satisfying \eqref{cond.order}. The following theorem of compactness is the key of this paper.
\begin{theorem}\label{thm.ScaleReduction}
	Let $A(x,y_1,\cdots,y_n) \in \mathcal{M}(\alpha,\beta;\Omega)$ satisfy \eqref{periodicity} and \eqref{smoothness}. For any sequence $\{ (\e^{(k)}_i )_{1\le i\le n}: k=1,2,\cdots \}$ satisfying \eqref{cond.order} with $\e^{(k)}_n \to 0$ as $k\to \infty$, there exists a subsequence $\{ (\e_i^{(k_j)})_{1\le i\le n}: j=1,2,\cdots \}$ and another coefficient matrix $B(x, y_1,y_2,\cdots, y_{m}) \in \mathcal{M}(\alpha,\beta;\Omega)$, satisfying \eqref{periodicity} and \eqref{smoothness}, with $1\le m\le n-1$ and $\delta = (\delta_1, \delta_2,\cdots, \delta_{m}) \in (0,\infty)^m$ such that
	\begin{equation*}
		A(x, x/\e^{(k_j)}_{1},\cdots,x/\e^{(k_j)}_{n}) \quad H\txt{-converges to}\quad B(x, x/\delta_1,\cdots, x/\delta_m),
	\end{equation*}
	as $j\to \infty$.
	
\end{theorem}

The point of the above theorem is not only the $H$-convergence (whose existence has been known \cite[Theorem 6.5]{Tbook09}), but also the fact that the limit is a matrix with the same periodic structure and strictly less number of variables. In Section 3, we will use this result to prove Theorem \ref{thm.main1} inductively.

To show the above theorem, we need the following lemma.
\begin{lemma}\label{lem.AjBj}
	Let $(A_j)_{j\ge 1}, (B_j)_{j\ge 1} \subset \mathcal{M}(\alpha,\beta;\Omega)$ be sequences of coefficient matrices. If $A_j$ $H$-converges to $\bar{A}$ and $\|B_j - A_j\|_{L^\infty(\Omega)} \to 0$ as $j\to \infty$, then $B_j$ $H$-converges to $\bar{A}$.
\end{lemma}
\begin{proof}
	Let $\Omega$ be a bounded Lipschitz domain and $\{ u_j \}$ be a sequence of weak solutions of $-\txt{div} (B_j \nabla u_j) = f$, bounded in $H_0^1(\Omega)$. Suppose $\bar{u}$ is the weak limit of $u_j$ in $H_0^1(\Omega)$ as $j\to \infty$. It suffices to show that $\bar{u}$ satisfies $-\txt{div} ( \bar{A} \nabla \bar{u})  = f$ in $\Omega$.
	
	First, observe that $u_j$ satisfies
	\begin{equation*}
		-\txt{div}  (A_j \nabla u_j) = f - \txt{div} ((A_j - B_j) \nabla u_j ) \qquad {\rm in } \,\,\Omega.
	\end{equation*}
	Write $u_j = v_j + w_j$ with $w_j$ being the solution of
	\begin{equation}\label{eq.vj}
		\left\{
		\begin{aligned}
			-\txt{div}  (A_j \nabla w_j) &= -\txt{div} ((A_j - B_j) \nabla u_j ) \qquad &\txt{in } \,\,\Omega,\\
			w_j &=0  \qquad &\txt{on } \partial\Omega.
		\end{aligned}	
		\right.
	\end{equation}
	Since $\{ u_j \}$ is bounded in $H^1(\Omega)$, by the assumption, $\norm{(A_j -B_j)\nabla u_j}_{L^2(\Omega)} \to 0$ as $j\to \infty$. The energy estimate for the Dirichlet problem (\ref{eq.vj}) implies $\norm{w_j}_{H^1(\Omega)} \to 0$ as $j\to \infty$. On the other hand, observe that $\{ v_j \}$ is a sequence of weak solutions (bounded in $H^1(\Omega)$) of $-\txt{div}  (A_j \nabla v_j) = f$. By the $H$-convergence of $A_j$, there is a subsequence $v_{j_k}$ so that $v_{j_k}$ converges to $\bar{v}$ weakly in $H^1$ and $\bar{v}$ is a weak solution of $-\txt{div} (\bar{A} \nabla \bar{v}) = f$ in $\Omega$. Finally, recall that $u_{j_k} = v_{j_k} + w_{j_k}$. By taking $k\to \infty$, we obtain $\bar{u} = \bar{v}$, which yields $-\txt{div}(\bar{A} \nabla \bar{u})  = f$ in $\Omega$. The proof is complete.
\end{proof}

The following is the reiterated homogenization theorem for quasi-periodic operators (see Section 5 for definition). It is important in our proof and of independent interest. The proof is nontrivial and will be postponed to Section 5.

\begin{theorem}\label{thm.quasiper.HC}
	If $A(x,y_1,y_2,\cdots,y_n) \in \mathcal{M}(\alpha,\beta;\Omega)$ is a $C^\gamma$-H\"{o}lder continuous (in all variables) coefficient matrix, quasi-periodic in every $y_k$ and $(\e_i)_{1\le i\le n}$ satisfies \eqref{cond.separate}, then $A_\e(x) = A(x,x/\e_1,x/\e_2,\cdots, x/\e_n)$ $H$-converges to some $C^\gamma$-H\"{o}lder continuous $\overline{A}(x) \in \mathcal{M}(\alpha,\beta;\Omega)$.
\end{theorem}

Now, we are ready to prove Theorem \ref{thm.ScaleReduction}.

\begin{proof}[Proof of Theorem \ref{thm.ScaleReduction}]
	In the following proof, for convenience, we will repeatedly take subsequences (sometime without mentioning) without changing the subscripts. By taking subsequences, we may assume $\e^{(k)}_i$ converges as $k\to \infty$ for each $i$. For convenience, we will denote the sequence $\{ (\e^{(k)}_i)_{1\le i\le n}: k =1,2,\cdots \}$ simply by $\{ \e^{(k)}_i\}$.
	
	Since $\{ \e^{(k)}_i\}$ is in a decreasing order in $i$, there are two integers $n'$ and $n''$ so that the following three cases take place as $k\to \infty$
	\begin{itemize}
		\item For $1\le i \le n'$, $\e^{(k)}_i \to \infty$;
		
		\item For $n'< i \le n''$, $\e^{(k)}_i \to \bar{\e}_{i} \in (0,\infty)$;
		
		\item For $n'' < i \le n$, $\e^{(k)}_i \to 0$.
	\end{itemize}
	Since we have assumed $\e^{(k)}_n \to 0$, then $n'' \le n-1$. But it is possible that $n' = 0$ or $n'' = 0$. A crucial fact about the first two cases is that since $\e^{(k)}_i$ converges to some nonzero limit (including infinity), the coefficient matrix hence is uniformly continuous with respect to the first $n''$ variables. Moreover, we have a stronger convergence (uniform convergence) for the first $n''$ variables, instead of homogenization ($H$-convergence), namely
	\begin{equation*}
		\bigg|A(x, \frac{x}{\e^{(k)}_{1}},\cdots,\frac{x}{\e^{(k)}_n}) - A(x, 0,\cdots, 0,\frac{x}{\bar{\e}_{n'+1}},\cdots, \frac{x}{\bar{\e}_{n''}}, \frac{x}{\e^{(k)}_{{n''+1}}},\cdots, \frac{x}{\e^{(k)}_n } ) \bigg| \to 0.
	\end{equation*}
	
	In view of Lemma \ref{lem.AjBj}, it suffices to study the $H$-convergence of
	\begin{equation*}
		\tilde{A}(x,\frac{x}{\bar{\e}_{n'+1}},\cdots, \frac{x}{\bar{\e}_{n''}},\frac{x}{\e^{(k)}_{{n''+1}}},\cdots, \frac{x}{\e^{(k)}_n }) := A(x, 0,\cdots, 0, \frac{x}{\bar{\e}_{n'+1}},\cdots, \frac{x}{\bar{\e}_{n''}}, \frac{x}{\e^{(k)}_{{n''+1}}},\cdots, \frac{x}{\e^{(k)}_n }).
	\end{equation*}
	Note that $\tilde{A}$ is uniformly $C^\gamma$-H\"{o}lder continuous with respect to the first $n''+1$ fixed variables and is 1-periodic with respect to the remaining $s:=n -n''$ faster variables. 	
	Recall that $\e^{(k)}_{{n''+i}} \to 0$ as $k\to \infty$ for each $i=1,2,\cdots, s$.
	To simplify our notation in the proof, we relabel $\e_{n''+i}^{(k)}$ as $\e^{(k)}_i$ (or just imagine $n'' = 0$). Then, we may rewrite $\tilde{A}$ as
	\begin{equation}\label{eq.tAs}
		\tilde{A}(x, \frac{x}{\bar{\e}_{n'+1}},\cdots, \frac{x}{\bar{\e}_{n''}}, \frac{x}{\e_{1}^{(k)}},\cdots, \frac{x}{\e_{s}^{(k)} } ),
	\end{equation}
	where $\e^{(k)}_i \ge \e_{i+1}^{(k)}$ and $\e^{(k)}_i \to 0$ as $k\to \infty$ for each $i$.
	
	To initiate our construction procedure, we first determine whether the scales of $\{ \e^{(k)}_i \} = \{ (\e^{(k)}_i)_{1\le i\le s} \}$ can be separated. By taking subsequence, we assume $\lim_{k\to \infty} \e_{i+1}^{(k)}/\e^{(k)}_i$ exists for each $i = 1,2,\cdots, s-1$. Since $\e^{(k)}_i \ge \e_{i+1}^{(k)}$, those limits are in $[0,1]$. We say the sequence $\{ \e^{(k)}_i \}$ and $\{ \e_{i+1}^{(k)} \}$ are not scale-separated (or they have the same scale), if $\lim_{k\to \infty} \e_{i+1}^{(k)}/\e^{(k)}_i >0$; otherwise, we say they are scale-separated. If all the sequence $\{ \e^{(k)}_i \}$ are scale-separated, then Theorem \ref{thm.quasiper.HC} directly implies $\tilde{A}(x,\frac{x}{\bar{\e}_{n'+1}},\cdots, \frac{x}{\bar{\e}_{n''}}, \frac{x}{\e_{1}^{(k)} },\cdots, \frac{x}{\e_{s}^{(k)} } )$ $H$-converges to some $\tilde{B}(x,\frac{x}{\bar{\e}_{n'+1}},\cdots, \frac{x}{\bar{\e}_{n''}})$. The difficulty arises when some sequences are not scale-separated. Our strategy to handle this difficulty is to reseparate the scales by repeated rewriting and rearrangement.
	
	Let $i_0$ be the largest integer so that $\{ \e_{i_0}^{(k)} \}$ and $\{ \e_{i_0-1}^{(k)} \}$ are not scale-separated. In other words, $\lim_{k\to \infty} \e_{i_0}^{(k)}/\e_{i_0-1}^{(k)}  >0 $, and for any $i > i_0$, $\{ \e^{(k)}_i \}$ and $\{ \e_{i-1}^{(k)} \}$ are scale-separated, i.e., $\lim_{k\to \infty} \e^{(k)}_i/\e_{i-1}^{(k)} = 0$. Let $S^+ = \{ \{\e_{i_0+1}^{(k)} \},\cdots, \{\e_{s}^{(k)} \} \}$. Let $S = \{  \{\e_{i_0-\ell+1}^{(k)}\},\cdots, \{\e_{i_0}^{(k)} \} \}$ be the largest set of sequences with the same scales (there are exactly $\ell$ sequences belonging to the same scales).
	Let $S^{-} = \{ \{\e^{(k)}_1 \},\cdots, \{\e_{i_0-\ell}^{(k)} \} \}$. Observe that $S^+$ contains sequences that all have separated scales, which are also separated from the sequences in $S$. While $S^{-}$ contains all the sequences (with the same or separated scales) with scales larger than the sequences in $S$. In the following, we will describe a procedure that reconstructs $S^{-}$ and $S$ so that at least one of the sequences  in $S$ could be moved to $S^+$. The inductive proof is complete if all the reconstructed sequences are moved to $S^+$. Note that initially, $S^+$ could be empty.
	
	Consider the sequences in $S$. By taking subsequences, assume $0<\gamma_{i} = \lim_{k\to \infty} \e_{i_0}^{(k)}/\e^{(k)}_i$ for $i_0 - \ell + 1 \le i \le i_0-1$. Thus, for each $i_0 - \ell + 1 \le i \le i_0-1$, we may break one scale into two scales by writing
	\begin{equation}\label{eq.x/e}
		\frac{x}{\e^{(k)}_i} = \Big( \frac{\e_{i_0}^{(k)}}{\e^{(k)}_i } - \frac{1}{\gamma_{i}} \Big) \frac{x}{\e_{i_0}^{(k)}} + \frac{1}{\gamma_{i}} \frac{x}{\e_{i_0}^{(k)}}.
	\end{equation}
	Note that $\{ \gamma_i \e_{i_0}^{(k)} \}$ is the same scale as $\{ \e_{i_0}^{(k)} \}$ and therefore the new expression \eqref{eq.x/e} does not increase the number of scales to $S$.
	Define $\tau_{i}^{(k)}:=\e_{i_0}^{(k)} /{\e^{(k)}_i } - 1/ \gamma_{i}$. Note that $\tau_{i}^{(k)} \to 0$ as $k\to \infty$. If $\{ \tau_k^{(k)} \}$ has a subsequence with all the  terms being zero, then ${\e^{(k)}_i} =\gamma_{i} \e_{i_0}^{(k)} $ in a subsequence. This turns out to be a simple case as we can reduce the number of scales in $S$ by combining $\e^{(k)}_i$ with $\e_{i_0}^{(k)}$. If all $\e_{i}^{(k)}$ with $i_0 - \ell + 1 \le i \le i_0-1$ can be reduced to $\e_{i_0}^{(k)}$, then $S$ only has one scale and $\tilde{A}$ is quasi-periodic at this scale. In this case, we may simply move this scale to $S^+$ as it is already separated from all the scales in $S^-$ and $S^+$, and then create a new set $S$ from $S^-$.
	
	Now, without loss of generality, we assume $\tau_{i}^{(k)} > 0$ (up to a subsequence) for all $i_0 - \ell + 1 \le i \le i_0-1$.
	Since $\tau_{i}^{(k)} \to 0$ as $k\to \infty$, $\tilde{\e}_{i}^{(k)}: = \e_{i_0}^{(k)}/\tau_{i}^{(k)}$ has a scale larger than $\{ \e_{i_0}^{(k)} \}$, namely, $\limsup_{k\to \infty} \e_{i_0}^{(k)}/\tilde{\e}_{i}^{(k)} = 0$. This observation is the key of our strategy for separating scales. Note that by (\ref{eq.x/e}),
	\begin{equation*}
		\begin{aligned}
			& \tilde{A}(*, \frac{x}{\e^{(k)}_{i_0-\ell+1}}, \cdots,\frac{x}{\e^{(k)}_{i_0-1}},\frac{x}{\e_{i_0}^{(k)}}, *) \\
			& = \tilde{A}(*, , \frac{x}{\tilde{\e}^{(k)}_{i_0-\ell+1}} +\frac{x}{\gamma_{i_0-\ell+1} \e_{i_0}^{(k)}}, \cdots, \frac{x}{\tilde{\e}^{(k)}_{i_0-1}} +\frac{x}{\gamma_{i_0-1} \e_{i_0}^{(k)}}, \frac{x}{\e_{i_0}^{(k)}}, *).
		\end{aligned}
	\end{equation*}
	For simplicity, we have used $*$ to represent the irrelevant variables in the above construction.
	
	Next, if we rewrite $\tilde{A}$ as
	\begin{equation*}
		\tilde{A}^{(1)}(*,y_{i_0-\ell+1},\cdots, y_{i_0-1}, y_{i_0},* ) = \tilde{A}(*, y_{i_0-\ell+1} + \frac{y_{i_0}}{\gamma_{i_0-\ell+1}},\cdots, y_{i_0-1} + \frac{y_{i_0}}{\gamma_{i_0-1}}, y_{i_0}, *),
	\end{equation*}
	then,
	\begin{equation*}
		\tilde{A}^{(1)}(*,\frac{x}{\tilde{\e}^{(k)}_{i_0-\ell+1}},\cdots, \frac{x}{\tilde{\e}^{(k)}_{i_0-1}}, \frac{x}{\e_{i_0}^{(k)}}, *) = \tilde{A}(*, \frac{x}{\e^{(k)}_{i_0-\ell+1}}, \cdots,\frac{x}{\e_{i_0}^{(k)}}, *).
	\end{equation*}
	Clearly, $\tilde{A}^{(1)}$ is still periodic with respect to $y_i$ with $i_0-\ell+1\le i\le i_0-1$. We also want to show that $\tilde{A}^{(1)}$ is quasi-periodic in $y_{i_0}$; see the definition in Section \ref{sec.quasipH}. Actually, if $\gamma_{i}/\gamma_{j}$ are all irrational for $i\neq j$, then $\tilde{A}^{(1)}(*,y_{i_0},*) = \tilde{A}(*,M_{i_0} y_{i_0},*)$ with  $M_{i_0} = (\gamma_{i_0-\ell+1}^{-1} I, \cdots, \gamma_{i_0-1}^{-1}  I, I)^T$. Here $I$ is the $d\times d$ identity matrix. In this case, it is not difficult to see that $M_{i_0}^T z \neq 0$ for all $z\in \Z^{\ell d}\setminus \{0\}$. If $\gamma_i/\gamma_j$ is rational for some $i\neq j$, then we can view $y_{i_0}/\gamma_i$ and $y_{i_0}/\gamma_j$ as one periodic variable (the period is a fixed constant, but may be large). With this idea, all $y_{i_0}/\gamma_j$ could be bundled into several periodic variables and the ratios of their periods are all irrational. Consequently, $\tilde{A}^{(1)}$ is quasi-periodic in $y_{i_0}$.
	
	Moreover, note that the sequence $\{ \e_{i_0}^{(k)} \}$ has a scale separated from all $\{ \tilde{ \e}_{i_0-\ell+1}^{(k)}\}, \cdots, \{\tilde{\e}_{i_0-1}^{(k)}\}$, as well as the sequences in $S^{-}$. Hence, we may now move $\{ \e_{i_0}^{(k)}\}$ to the set $S^+$. As a result, the number of scales in $S$ is reduced at least by one.
	
	Finally, to repeat the above process, we need to reconstruct $S$ and $S^{-}$. Note that in the above construction, the sequence $\{ \tilde{\e}_{i}^{(k)} \}$ with $i_0-\ell+1 \le i \le i_0 - 1$ may not converge to zero. As previously, the sequences converging to positive constants or infinity can be viewed as slow variables and handled by Lemma \ref{lem.AjBj}. The remaining sequences of ${ \tilde{\e}^{(k)}_{i} }$ (which converge to zero), together with the sequences in $S^{-}$ form a new set of sequences (at most $i_0-1$ sequences), from which we can reconstruct new sets $S^-$ and $S$. Thus, we can repeat the previous process at most $i_0-1$ times and finally move all the sequences into $S^+$ so that all the scales are separated. Consequently, there exists a quasi-periodic matrix $\tilde{A}^{\sharp}$ so that (up to a subsequence)
	\begin{equation}\label{eq.As-A}
		\tilde{A}^{\sharp}(x,\frac{x}{\delta_1}, \frac{x}{\delta_2},\cdots, \frac{x}{\delta_m}, \frac{x}{\hat{\e}^{(k)}_{m+1}},\cdots, \frac{x}{\hat{\e}^{(k)}_{m+m_0}}) = A(x,\frac{x}{\e^{(k)}_{1}},\cdots, \frac{x}{\e^{(k)}_n}) + o(1)
	\end{equation}
	as $k\to \infty$. We point out that $\delta_i$ are the positive limits $\bar{\e}_i >0$ or the positive limits of $\{ \tilde{\e}^{(k)}_{i} \}$ in the intermediate steps. Also, by the above construction, $\tilde{A}^\sharp$ is 1-periodic and $C^\gamma$-H\"{o}lder continuous with respect to the first $m$ oscillating variables. Moreover, the sequence $\{ \hat{\e}^{(k)}_{m+i}\}$ with $1\le i\le m_0$ are all convergent to zero and scale-separated. Since we assume $\e^{(k)}_n \to 0$, the final set $S^+$ cannot be empty. Thus, $m_0 \ge 1$ and therefore $m+m_0 \le n$ implies $m \le n-1$.
	
	Finally, by Theorem \ref{thm.quasiper.HC},
	\begin{equation}\label{eq.Asharp.H}
		\tilde{A}^{\sharp}(x, \frac{x}{\delta_1}, \frac{x}{\delta_2},\cdots, \frac{x}{\delta_m}, \frac{x}{\hat{\e}^{(k)}_{m+1}},\cdots, \frac{x}{\hat{\e}^{(k)}_{m+m_0}}) \quad H\txt{-converges to some } B(x, \frac{x}{\delta_1}, \frac{x}{\delta_2},\cdots, \frac{x}{\delta_m}).
	\end{equation}
	Hence, the desired $H$-convergence follows from (\ref{eq.Asharp.H}), (\ref{eq.As-A}) and Lemma \ref{lem.AjBj}. To see the periodic structure  of $B(x,y_1,\cdots, y_m)$, it is sufficient to observe the integer translation invariance with respect to $x/\delta_i$ in \eqref{eq.Asharp.H} (for more details, see the construction of $B$ in proof of Theorem \ref{quasi-th4}). The condition \eqref{smoothness} also follows easily by noting that the regularity is preserved under the $H$-convergence $\eqref{eq.Asharp.H}$.
\end{proof}

\begin{remark}\label{rmk.cpt}
	In the above argument, the fixed matrix $A(x, y_1,\cdots, y_n)$ can   be replaced by a sequence $ \{ A^{(k)}(x, y_1,\cdots, y_n) \}$ as long as they satisfy the same assumptions \eqref{ellipticity} - \eqref{smoothness} with the same constants. This is because the set of uniformly bounded and equicontinuous functions on $\Omega\times \T^{d\times n}$ is compact under uniform convergence (Arzel\`{a}-Ascoli Theorem). Thus, there is a subsequence of $\{ A^{(k)}(x, y_1,\cdots, y_n) \}$ converging uniformly to some fixed $A(x,y_1,\cdots, y_n)$ and Lemma \ref{lem.AjBj} applies.
\end{remark}

\section{Stable $C^\alpha$ estimate}
In this section, we use an upgraded method of compactness to prove Theorem \ref{thm.main1} inductively on the number of scales $n$.
\begin{proof}[Proof of Theorem \ref{thm.main1}]
	By normalization, assume $(\fint_{B_1}|u_\e|^2)^{1/2} = 1$. Let $j$ be the number of scales of the coefficient matrix $A_\e = A(x, x/\e_1,\cdots, x/\e_j)$. We prove by induction on $j$. For $j=1$, $A_\e(x) = A(x, x/\e_1)$ is a 1-scale periodic coefficient matrix. Then the uniform regularity (\ref{est.mainCa}) for any $\alpha\in (0,1)$ is well-known \cite{al87,nsxjfa2020}. Now, suppose that (\ref{est.mainCa}) holds for all $1\le j\le n-1$ with constant $C_{n-1,\alpha}$.
	In other words, if $u_\e$ is a normalized solution of
	(\ref{eq.main}), then for any $\alpha\in (0,1)$, there exists $C_{n-1,\alpha}>0$ such that
	\begin{equation}\label{est.j<=n-1}
		\inf_{q\in \R} \bigg( \fint_{B_\theta} |u_{\e} - q|^2 \bigg)^{1/2} \le C_{n-1,\alpha} \theta^\alpha,
	\end{equation}
	for all $\theta\in (0,1)$. Note that (\ref{est.j<=n-1}) is a straightforward consequence of (\ref{est.mainCa}).
	
	 We would like to show (\ref{est.mainCa}) for $j = n$. Recall that $(\e_1,\e_2,\cdots, \e_n)$ satisfies $\e_k \ge \e_{k+1}$, $1\leq k\leq n-1$.
	
	\textbf{Step 1: One-step improvement.}
	
	We claim that there exist $\delta_0,\theta>0$ such that if $\e_n < \delta_0$, then
	\begin{equation}
	\inf_{q\in \R} \bigg( \fint_{B_\theta} |u_{\e} - q|^2 \bigg)^{1/2} \le \theta^\alpha.
	\end{equation}
	If the above claim is not true, then for any given $\theta>0$ (to be determined), there exist a sequence $\{(\e^{(k)}_1, \cdots, \e^{(k)}_{n} ) \}$, with $\e^{(k)}_{n} $ converging to zero, a sequence of normalized solutions  $\{ u_{\e^{(k)}}\}$ to (\ref{eq.main}) corresponding to $A_{\e^{(k)}}$, so that for all $k \ge 1$
	\begin{equation}\label{cond.contra}
	\inf_{q\in \R} \bigg( \fint_{B_\theta} |u_{\e^{(k)}} - q|^2 \bigg)^{1/2} > \theta^\alpha.
	\end{equation}
	Since $(\fint_{B_1}|u_{\e^{(k)} }|^2)^{1/2}  = 1$, the Caccioppoli inequality implies that $\norm{ u_{\e^{(k)}}}_{H^1(B_{3/4})} \le C$. Hence by compactness, we can find a subsequence, still denoted by $\{ u_{\e^{(k)}} \}$, converging to $u_0$ strongly in $L^2(B_{3/4})$, i.e.,
	\begin{equation}\label{est.L2conv}
		\norm{u_{\e^{(k)}}  - u_0}_{L^2(B_{3/4})} \to 0, \qquad \txt{as } k\to \infty.
	\end{equation}
	By Theorem \ref{thm.ScaleReduction}, there exists a subsequence of $\{(\e^{(k)}_1, \cdots, \e^{(k)}_{n} ) \}$, still denoted by $\{(\e^{(k)}_1, \cdots, \e^{(k)}_{n} ) \}$, so that $A_{\e^{(k)}}$ $H$-converges to $B_\delta(x) = B(x, x/\delta_1,x/\delta_2,\cdots,x/\delta_{m})$ where $1\le m\le n-1$ and $B(x, y_1,y_2,\cdots,y_m)$ is a 1-periodic coefficient matrix in $y_i$. Therefore, $u_0$ is a weak solution of (because of $H$-convergence)
	\begin{equation}\label{eq.Bdelta.u0}
	-\txt{div} (B_{\delta} \nabla u_0) = 0, \qquad \txt{ in } B_{3/4}.
	\end{equation}
	Note that even though the $H$-convergence is defined through Dirichlet problems in Definition \ref{H-convergence}, it is actually not related to any particular boundary conditions; see \cite[Lemma 10.3]{Tbook09}. This is why \eqref{eq.Bdelta.u0} holds, regardless of the boundary conditions of $u_{\e^{(k)}}$ on $\partial B_{3/4}$.

	Now, since $B_\delta$ has at most $n-1$ scales, by our inductive assumption and (\ref{est.j<=n-1}), we know that for any $\theta \in (0,1/2)$,
	\begin{equation*}
	\inf_{q\in \R} \bigg( \fint_{B_\theta} |u_0 - q|^2 \bigg)^{1/2} \le C_{n-1,\beta} \theta^\beta.
	\end{equation*}
	Choose $\beta \in (\alpha,1)$ and $\theta$ small so that $C_{n-1,\beta} \theta^\beta < \theta^\alpha$. For such $\theta$ (fixed), (\ref{est.L2conv}) and the above estimate imply
	\begin{equation*}
	\lim_{k\to \infty} \inf_{q\in \R} \bigg( \fint_{B_\theta} |u_{\e^{(k)}} - q|^2 \bigg)^{1/2} = \inf_{q\in \R} \bigg( \fint_{B_\theta} |u_0 - q|^2 \bigg)^{1/2} < \theta^\alpha.
	\end{equation*}
	This contradicts to (\ref{cond.contra}). Thus, we have proved the claim.
	
	\textbf{Step 2: Iteration.} Fix $\delta_0>0$ and $\theta>0$ obtained from Step 1. The scale-invariant estimate in Step 1 can be written as
	\begin{equation*}
		\inf_{q\in \R} \bigg( \fint_{B_\theta} |u_{\e} - q|^2 \bigg)^{1/2} \le \theta^\alpha \inf_{q\in \R} \bigg( \fint_{B_1} |u_{\e} - q|^2 \bigg)^{1/2}.
	\end{equation*}
	We added an arbitrary $q$ on the right-hand side because $u_\e -q$ is also a weak solution.
	Now, assume for $s>0$, $u_{\e,s}(x) = u_\e(x/s)$. By rescaling, $u_{\e,\theta^{-k}}$ is a solution for (\ref{eq.main}) with the coefficient matrix $A_{\e}(\theta^{k}\cdot)$. Thus, if $\theta^{-k} \e_n < \delta_0$,
	\begin{equation*}
		\inf_{q\in \R} \bigg( \fint_{B_\theta} |u_{\e,\theta^{-k}} - q|^2 \bigg)^{1/2} \le \theta^\alpha \inf_{q\in \R} \bigg( \fint_{B_1} |u_{\e,\theta^{-k} } - q|^2 \bigg)^{1/2}.
	\end{equation*}
	By rescaling,
	\begin{equation*}
		\inf_{q\in \R} \bigg( \fint_{B_{\theta^{k+1}}} |u_{\e} - q|^2 \bigg)^{1/2} \le \theta^\alpha \inf_{q\in \R} \bigg( \fint_{B_{\theta^k}} |u_{\e} - q|^2 \bigg)^{1/2}.
	\end{equation*}
	Consequently, if $\theta^{-k } \e_n < \delta_0$
	\begin{equation*}
		\inf_{q\in \R} \bigg( \fint_{B_{\theta^{k+1}}} |u_{\e} - q|^2 \bigg)^{1/2} \le \theta^{(k+1)\alpha} \inf_{q\in \R} \bigg( \fint_{B_1} |u_{\e} - q|^2 \bigg)^{1/2}.
	\end{equation*}
	This implies that for all $\e_n < \delta_0$ and $r > \e_n/\delta_0$,
	\begin{equation}\label{est.Br-B1}
		\inf_{q\in \R} \bigg( \fint_{B_{r}} |u_{\e} - q|^2 \bigg)^{1/2} \le Cr^\alpha \bigg( \fint_{B_1} |u_{\e}|^2 \bigg)^{1/2}.
	\end{equation}
	\textbf{Step 3: Blow-up.} Note that $\e_k \ge \e_n$ for all $k < n$. If $\e_n > \delta_0$, then $A_\e$ is uniformly H\"{o}lder continuous, which implies the desired estimate by the classical Schauder estimate without the previous two steps. If $\e_n < \delta_0$ and $r< \e_n/\delta_0$, then we may apply the classical Schauder estimate to the rescaled solution $u_{ \e, \delta_0/\e_n}$, together with \eqref{est.Br-B1}. Precisely, we have
	\begin{equation}\label{est.Br-B1-2}
		\begin{aligned}
		\inf_{q\in \R} \bigg( \fint_{B_{r}} |u_{\e} - q|^2 \bigg)^{1/2} &= \inf_{q\in \R} \bigg( \fint_{B_{r\delta_0/\e_n}} |u_{\e ,\delta_0/\e_n } - q|^2 \bigg)^{1/2} \\
		& \le C(r\delta_0/\e_n)^\alpha \inf_{q\in \R} \bigg( \fint_{B_1} |u_{\e,\delta_0/\e_n } - q|^2 \bigg)^{1/2} \\
		& = C(r\delta_0/\e_n)^\alpha \inf_{q\in \R} \bigg( \fint_{B_{\e_n/\delta_0}} |u_{\e} - q|^2 \bigg)^{1/2} \\
		& \le Cr^\alpha \bigg( \fint_{B_1} |u_{\e}|^2 \bigg)^{1/2},
	\end{aligned}
\end{equation}
where we have used (\ref{est.Br-B1}) with $r = \e_n/\delta_0$ in the last inequality. Hence, (\ref{est.Br-B1}) in fact holds for all $r\in (0,1)$. Finally, since this estimate is translation invariant, we may show (\ref{est.Br-B1}) in any balls centered at $x\in B_{1/2}$. This implies (\ref{est.mainCa}) by a standard result of Campanato spaces.
\end{proof}

\begin{remark}
	In view of Remark \ref{rmk.cpt}, the constant $C$ in Theorem \ref{thm.main1} depends at most on $d, n, \Lambda, L$ and $\gamma$.
\end{remark}

\section{Homogenization under  variable-separation condition}

 This part is devoted to quantitative homogenization of $ -\txt{div} ( A_\va(x) \na )$ under the variable-separation condition. Precisely speaking, we consider the special case
\begin{align}
	A_\va(x)  = A(x, x_1/\va_1, x_2/\va_2,\cdots ,x_d/\va_d),
\end{align}
where $A=A(x,y), y=(y_1,\dots ,y_d) \in \R^d$. Without loss of generality, we assume $1\ge \e_1 \ge \e_2 \ge \cdots \ge \e_d > 0$. Suppose that $A$
  satisfies the ellipticity condition \eqref{ellipticity}, periodicity condition \eqref{periodicity} in $y$ and smoothness condition \eqref{smoothness} in both $x$ and $y$. We will show that not only the solution $u_\va$ of $-\text{div} ( A_\e(x)  \na u_\va )=0$ admits the uniform Lipschitz estimate, but also it can be approximated in $L^2(\Omega)$ with error $O(\e_1)$ by the solution of an elliptic equation with non-oscillating coefficients.

\subsection{Reperiodization and stable Lipschitz estimates} \label{sec.reperiodization}
We now introduce a method of reperiodization to establish the uniform Lipschitz estimate.
The method is based on a simple fact: if $A(x,y_1, y_2, \cdots, y_d)$ is a 1-periodic function in $y_i$'s, then $A(x,m_1 y_1, m_2 y_2, \cdots, m_d y_d)$ is also a 1-periodic function in $y_i$'s, provided that $m_i$'s are all integers. Moreover, even if $m_i \ge 1$ are not integers, we can still reduce it to a 1-periodic function by a non-degenerate change of variables.

To prove Theorem \ref{thm.main2}, we need the following result as a black box. Its proof may be found in \cite{nsxjfa2020}.
\begin{theorem}\label{thm.2Scale.Lip}
	Let  $u_\e \in H^1(B_1)$ be a weak solution of $-\txt{div} (A(x,x/\e) \nabla u_\e) = 0$ in $B_1$, with $A(x,y)$ satisfying conditions  \eqref{ellipticity}, \eqref{periodicity} (with $n=1$), and
	\begin{equation}\label{ccon}
		|A(x,y) - A(x',y)| \le L|x-x'|^\gamma.
	\end{equation} Then
	\begin{equation*}
		\sup_{\e \le r \le 1} \bigg( \fint_{B_r} |\nabla u_\e|^2 \bigg)^{1/2} \le C\bigg( \fint_{B_1} |\nabla u_\e|^2 \bigg)^{1/2},
	\end{equation*}
where $C$ depends only on $d, \Lambda, L$ and $ \gamma$.
\end{theorem}

It is crucial that in the above result there is no smoothness assumption for $A(x,y)$ in $y$, for which we only need the periodicity structure.
Theorem \ref{thm.main2} will be proved by iteration. The first step is to prove the Lipschitz estimate at the largest scale $r\ge \e_1$. To this end, we rewrite $A_\e(x) = A(x, x_1/\va_1, x_2/\va_2,\cdots ,x_d/\va_d)$ in a form satisfying the assumptions in Theorem \ref{thm.2Scale.Lip}. Since $\e_1 = \max_{1\le i\le d} \e_i$, we may write
\begin{equation*}
	A_\e(x) = A\big(x, \frac{x_1}{\e_1}, \frac{\lambda_2 x_2}{ \e_1},\cdots, \frac{\lambda_d x_d}{\e_1} \big),
\end{equation*}
where $\lambda_i = \e_1/\e_i \ge 1$. Now, define
\begin{equation}\label{Mlm}
		M_\lm=
		\left(
		\begin{array}{cccccccc}
			\lambda_1 &  0 & 0  & \cdots & 0  \\
			0 & \lambda_2 & 0  & \cdots &0 \\
			0 & 0& \lambda_3 & \cdots& 0   \\
			\vdots & \vdots & \vdots & \ddots& \vdots \\
			0 & 0  & 0 &\cdots   & \lambda_d
		\end{array}
		\right).
\end{equation}
	We may write $A_\e(x) = A(x,M_\lambda x/\e_1)$. Notice that $A(x,M_\lambda y)$ is not $1$-periodic in $y$. Now, we would like to reperiodize $A$ by a change of variables. Let
	\begin{equation*}
		\fl{M_\lm}=
		\left(
		\begin{array}{cccccccc}
			\fl{\lambda_1} &  0 & 0  & \cdots & 0  \\
			0 & \fl{\lambda_2} & 0  & \cdots &0 \\
			0 & 0& \fl{\lambda_3} & \cdots& 0   \\
			\vdots & \vdots & \vdots & \ddots& \vdots \\
			0 & 0  & 0 &\cdots   & \fl{\lambda_d}
		\end{array}
		\right),
	\end{equation*}
where $\fl{t}$ denotes the integer part of $t \in \R$ (i.e., the largest integer not exceeding $t$). Note that
\begin{equation}\label{est.PhiPositivity}
	I_{d\times d}\le  \Phi_\lambda:=M_\lambda \fl{M_\lambda}^{-1} \le 2 I_{d\times d}.
\end{equation}

Consider
\begin{equation*}
	A^\sharp(x,y) =  \Phi_\lambda A(\Phi^{-1}_\lambda x, \fl{M_\lm}y) \Phi_\lambda.
\end{equation*}
The key observation is that $A^\sharp$ satisfies all the conditions of $A$ in Theorem \ref{thm.2Scale.Lip} if $A(x,y)$ does. In particular, $A^\sharp(x,y)$ is 1-periodic in $y$ because all the diagonal entries in $\fl{M_\lm}$ are integers. Notice that $A^\sharp(x,y)$ is still oscillating in $y$ since smaller periodic structures for variable $y_i$ with $2\le i\le d$ are compressed into a single 1-period. However, this is harmless since, as we have mentioned, Theorem \ref{thm.2Scale.Lip} requires no smoothness of $A^\sharp(x,y)$ in $y$.

Let $u_\e$ be a weak solution of $-\txt{div} (A_\e(x) \nabla u_\e) = 0$ in $B_1$. Then $v_\e(x) = u_\e(\Phi^{-1}_\lambda x)$ satisfies
\begin{equation*}
	-\txt{div} (A^\sharp_\e \nabla v_\e) = 0  \quad \text{in } \Phi_\lambda(B_1),
\end{equation*}
where $A^\sharp_\e = A^\sharp(x,x/\e_1)$. As a consequence of Theorem \ref{thm.2Scale.Lip}, we have
\begin{equation*}
	\sup_{ \e_1 \le r \le 1} \bigg( \fint_{B_r} |\nabla v_\e|^2 \bigg)^{1/2} \le C\bigg( \fint_{B_1} |\nabla v_\e|^2 \bigg)^{1/2}.
\end{equation*}
In view of \eqref{est.PhiPositivity}, this implies
\begin{lemma}\label{lem.e1}
 	Let $u_\e$ be a weak solution of $-\txt{div} (A_\e(x) \nabla u_\e) = 0$ in $B_1$, where $A=A(x,y)$ satisfies the assumptions \eqref{ellipticity}, \eqref{periodicity} and \eqref{ccon}. Then
	\begin{equation}
		\sup_{ \e_1 \le r \le 1} \bigg( \fint_{B_r} |\nabla u_\e|^2 \bigg)^{1/2} \le C\bigg( \fint_{B_1} |\nabla u_\e|^2 \bigg)^{1/2}.
	\end{equation}
\end{lemma}

\begin{proof}[Proof of Theorem \ref{thm.main2}]
	By translation and dilation, we only need to show $$|\nabla u_\e(0)| \le C\norm{ u_\e}_{L^2(B_1)}.$$ We claim that Lemma \ref{lem.e1} implies
	\begin{equation}\label{est.eiStep}
		\sup_{\e_{i+1} \le r\le \e_i } \bigg( \fint_{B_r} |\nabla u_\e|^2 \bigg)^{1/2} \le C\bigg( \fint_{B_{\e_i}} |\nabla u_\e|^2 \bigg)^{1/2}.
	\end{equation}
	Actually, this follows from Lemma \ref{lem.e1} and a blow-up argument. Precisely, let
	\begin{equation*}
	A^{i}_\e(x) = A_\e(\e_i x) = A(\e_i x, \mu_1 x_1, \cdots, \mu_{i-1}x_{i-1}, x_i, \frac{x_{i+1}}{\e_{i+1}/\e_i }, \cdots, \frac{x_d}{\e_d/\e_i  }),	
	\end{equation*}
	where $\mu_k = \e_i/\e_k$.
	Because $\e_i, \mu_1,\cdots, \mu_{i-1} \le 1$, the involved variables are all slow variables ($A^i_\e$ is uniformly $C^\gamma$-H\"{o}lder continuous in these variables) and we only have $d-i$ fast variables. Note that the largest scale now is $\e_{i+1}/\e_i$. If we put $w_\e(x) = u_\e(\e_i x)$, then $-\txt{div} (A^i_\e(x) \nabla w_\e) = 0$ in $B_1$. Thus, by Lemma \ref{lem.e1},
	\begin{equation*}
		\sup_{ \e_{i+1}/\e_i \le r \le 1} \bigg( \fint_{B_r} |\nabla w_\e|^2 \bigg)^{1/2} \le C\bigg( \fint_{B_1} |\nabla w_\e|^2 \bigg)^{1/2}.
	\end{equation*}
	By rescaling, this leads to (\ref{est.eiStep}).
	
	Now, by Lemma \ref{lem.e1} and iterating (\ref{est.eiStep}), we have
	\begin{equation}\label{est.ed}
		\sup_{ \e_d \le r \le 1/2} \bigg( \fint_{B_r} |\nabla u_\e|^2 \bigg)^{1/2} \le C\bigg( \fint_{B_{1/2}} |\nabla u_\e|^2 \bigg)^{1/2}.
	\end{equation}
 	Finally, another blow-up argument yields
	 \begin{equation}\label{est.0toed}
		|\nabla u_\e(0)| \le \bigg( \fint_{B_{\e_d}} |\nabla u_\e|^2 \bigg)^{1/2}.
	\end{equation}
	The desired estimate then follows directly from \eqref{est.ed}, \eqref{est.0toed} and the Caccioppoli inequality. 
\end{proof}

\subsection{Applications to parabolic equations}

The above approach of reperiodization can be applied to parabolic equations with coefficients depending on time
\begin{equation}\label{eq.parabolic}
	\frac{\partial u_\e}{\partial t} - \txt{div} (A_\e(x,t) \nabla u_\e) = 0 \quad \text{in } B_1\times (-1,0).
\end{equation}
Let us first consider a slightly simpler case
$A_\e(x,t) = A(x,t, x/\e, t/\de^2)$ with arbitrary $(\e,\delta)\in (0,1]^2$, and $A(x,t,y,s)$ is periodic in $(y,s) \in \R^d\times \R$. We recall that the special cases $A_\e(x,t) = A(x/\e, t/\e^k)$ with $ k>0$, and $A_\e(x,t) = A(x/\e,t/\delta^2)$ with arbitrary $(\e,\delta)\in (0,1]^2$ have been studied in \cite{GS20} and \cite{GN21}, respectively. In particular, the Lipschitz estimate uniform in $(\e,\delta)$ was derived in \cite{GN21} by the access decay method combined with a proper rescaling argument. 
In this subsection, we would like to show how the technique of reperiodization can be used to reduce the general case to the classical self-similar case $\delta = \e$ (one-scale), and then generalize the result to the spatially multiscale situation. In the following, we define $Q_r = Q_r(0,0) = B_r  \times (-r^2, 0)$.

First of all, as in elliptic case, we need the following basic result:
\begin{theorem}[\cite{GN21}]\label{thm.2Scale.Parabolic}
	Assume $A(x,t,y,s)$ satisfy 
the ellipticity condition: $\exists \Lambda\in (0,1]$ such that  
	\begin{equation}\label{pell.para}
		\Lambda |\xi|^2 \le A(x,t,y,s)\xi \cdot \xi \le \Lambda^{-1} |\xi|^2   
	\end{equation}for any  $(x,t), (y,s) \in\R^{d+1}$;
  the periodicity condition:     
	\begin{equation}\label{pcon.para}
		A(x,t,y+z,s+\tau) = A(x,t,y,s) 
	\end{equation}  for any  $(x,t) , (y,s) \in\R^{d+1}$ and any  $(z,\tau)\in \Z^{d+1}$;
	and the smoothness condition: $\exists L>0,$ and $\gamma \in (0,1]$ such that 
	\begin{equation}\label{ccon.para}
		|A(x,t,y,s) - A(\tilde{x},\tilde{t},y,s)| \le L\big(|x-\tilde{x}| + {|t-\tilde{t}|}^{\frac12} \big)^\gamma
	\end{equation} for any  $(x,t),(\tilde{x},\tilde{t}),$ and $ (y,s) \in\R^{d+1}$.
	Let  $u_\e \in L^2( (-1, 0); H^1(B_1))$ be a weak solution of \eqref{eq.parabolic} with $A_\e(x,t) = A(x, t, x/\e,t/\e^2)$. Then
	\begin{equation*}
		\sup_{\e \le r \le 1} \bigg( \fint_{Q_r} |\nabla u_\e|^2 \bigg)^{1/2} \le C\bigg( \fint_{Q_1} |\nabla u_\e|^2 \bigg)^{1/2},
	\end{equation*}
	where $C$ depends only on $d, \Lambda, L$ and $\gamma$.
\end{theorem}

Again, we point out that in the above theorem $A(x,t,y,s)$ has no smoothness assumption on $(y,s)$. Now, we are going to use the idea of reperiodization to reduce the general case $A_\e(x,t) = A(x,t,x/\e,t/\delta^2)$ to the special case contained in Theorem \ref{thm.2Scale.Parabolic}. This method works for parabolic equations essentially because $t$ and $x$ are independent variables satisfying the variable-separation condition. Due to the non-symmetry between the $x$ and $t$ variables, we have to consider two cases based on the ratio $\e/\delta$.

\textbf{Case 1:} Let us first assume $\lambda:=\delta/\e \ge 1$. Consider the coefficient matrix
\begin{equation}\label{scaled.co}
	A^\sharp(x,t,y,s) =\lambda^2 \fl{\lambda}^{-2} A(\fl{\lambda} \lambda^{-1} x,t,\fl{\lambda}y,s).
\end{equation}
Note that $1\le \lambda \fl{\lambda}^{-1} \le 2$. It is easy to see that $A^\sharp(x,t,y,s)$ satisfies the assumptions on $A$ in Theorem \ref{thm.2Scale.Parabolic}. Particularly, it is periodic in $(y,s)$. Moreover, by setting $v_\e(x,t) = u_\e (\fl{\lambda} \lambda^{-1} x,t)$, we see that $u_\e$ is a weak solution of \eqref{eq.parabolic} with $A_\e(x,t) = A(x,t, x/\e, t/\de^2)$ if and only if $v_\e$ satisfies
\begin{equation}
	\frac{\partial v_\e}{\partial t} - \txt{div} (A^\sharp(x,t,x/\delta,t/\delta^2) \nabla v_\e) = 0 \quad \text{in } (\lambda\fl{\lambda}^{-1} B_1)\times ( -1, 0).
\end{equation}
Therefore, Theorem \ref{thm.2Scale.Parabolic} implies
\begin{equation*}
	\sup_{\delta \le r \le 1} \bigg( \fint_{Q_r} |\nabla v_\e|^2 \bigg)^{1/2} \le C\bigg( \fint_{Q_1} |\nabla v_\e|^2 \bigg)^{1/2}.
\end{equation*}
By changing variables back to $u_\e$, we have
\begin{equation}\label{est.delta-1}
	\sup_{\delta \le r \le 1} \bigg( \fint_{Q_r} |\nabla u_\e|^2 \bigg)^{1/2} \le C\bigg( \fint_{Q_1} |\nabla u_\e|^2 \bigg)^{1/2}.
\end{equation}

\begin{remark}
To obtain the same estimate for $\e < r < \delta$, we need to further assume that $A(x,t,y,s)$ is H\"{o}lder continuous with respect to $s$. By rescaling $w_\e(x,t) = u_\e(\delta^{-1} x, \delta^{-2} t)$, the equation \eqref{eq.parabolic} can be written as
\begin{equation}\label{eq.parabolic-1}
	\frac{\partial w_\e}{\partial t} - \txt{div} (A(\delta x,\delta^2 t, \lambda x, t) \nabla w_\e) = 0 \quad \text{in } \delta^{-1}B_1\times \delta^{-2}( -1, 0).
\end{equation}
Now, note that $A(\delta x,\delta^2 t, \lambda x, t)$ is uniformly H\"{o}lder continuous with respect to the first, second and fourth variables and periodically oscillating only with respect the third variable since $\lambda \ge 1$. Consequently, we apply Theorem \ref{thm.2Scale.Parabolic} to \eqref{eq.parabolic-1} to obtain
\begin{equation*}
	\sup_{\lambda^{-1} \le r \le 1 } \bigg( \fint_{Q_r} |\nabla w_\e|^2 \bigg)^{1/2} \le C\bigg( \fint_{Q_1} |\nabla w_\e|^2 \bigg)^{1/2}.
\end{equation*}
Recall that $\lambda = \delta/\e$. By rescaling again, we obtain
\begin{equation*}
	\sup_{\e \le r \le \delta} \bigg( \fint_{Q_r} |\nabla u_\e|^2 \bigg)^{1/2} \le C\bigg( \fint_{Q_\delta} |\nabla u_\e|^2 \bigg)^{1/2}.
\end{equation*}
This, combined with \eqref{est.delta-1}, gives the large-scale Lipschitz estimate down to the scale $r \ge \e$. Finally, for $0<r<\e$, a blow-up argument applies if $A$ is H\"{o}lder continuous with respect to all the variables, and we may obtain $|\nabla u_\e(0,0)| \le C\norm{ \nabla u_\e}_{L^2(Q_1)}$. By translation invariance, this implies $\norm{ \nabla u_\e}_{L^\infty(Q_{1/2})} \le  C\norm{ \nabla u_\e}_{L^2(Q_1)}$.
\end{remark}

\textbf{Case 2:} Now assume $\e/\delta > 1$ and put $\lambda = \e^2/\delta^2$. The argument  for the Lipschitz estimate is similar to Case 1. In this case,  instead of \eqref{scaled.co} we consider the  coefficient matrix
\begin{equation*}
	A^\flat(x,t,y,s) := \fl{\lambda} \lambda^{-1}  A(x,\fl{\lambda} \lambda^{-1} t, y, \fl{\lambda} s).
\end{equation*}
This new matrix is 1-periodic in $(y,s)$ and by a change of variable $(t,x) \mapsto (\fl{\lambda} \lambda^{-1} t, x)$, the equation \eqref{eq.parabolic} can be written equivalently as
\begin{equation}
	\frac{\partial v_\e}{\partial t} - \txt{div} (A^\flat(x,t,x/\e,t/\e^2) \nabla v_\e) = 0 \quad \text{in }  B_1\times  (\lambda\fl{\lambda}^{-1}(-1,0)).
\end{equation}
By mimicking the argument in Case 1, we can obtain (\ref{est.delta-1}) with $\delta \le r\le 1$ replaced by $\e\le r \le 1$.

In view of the results derived in Case 1 and Case 2, we get 
\begin{equation}\label{est.Lip.e-d}
	\sup_{\delta\vee \e \le r \le 1} \bigg( \fint_{Q_r} |\nabla u_\e|^2 \bigg)^{1/2} \le C\bigg( \fint_{Q_1} |\nabla u_\e|^2 \bigg)^{1/2},
\end{equation}
where $a\vee b = \max \{ a,b\}$, proided that $A$ satisfies the conditions \eqref{pell.para}-\eqref{ccon.para}. It is important to notice that there is no smoothness assumption for $A(x,t,y,s)$ on $y$ or $s$.

Now, we consider the multiscale coefficients with the spatial  variable  separation, namely, $A_\e(x,t) = A(x,t,x_1/\e_1,\cdots,x_d/\e_d, t/\delta^2)$  with arbitrary $(\e_i)_{1\le i\le d} \in (0,1]^d$ and $\delta \in (0,1]$.  Without loss of generality, we assume $1 \ge \e_1 \ge \e_2 \ge \cdots \ge \e_d > 0$. We claim the following Lipschitz estimate.

\begin{theorem}\label{thm.parabolic}
	Assume $A(x,t,y,s)$ satisfies conditions \eqref{pell.para},   \eqref{pcon.para}, and 
	\begin{equation*}
		|A(x,t,y,s) - A(\tilde{x},\tilde{t}, \tilde{y},\tilde{s})| \le L\big(|x-\tilde{x}| + |y-\tilde{y}| + |t-\tilde{t}|^{\frac12} + |s-\tilde{s}|^{\frac12} \big)^{\gamma},
	\end{equation*}
	Let $u_\e \in L^2((-1, 0); H^1(B_1))$ be a weak solution of \eqref{eq.parabolic}, with
	\begin{equation*}
		A_\e(x,t) = A(x,t,x_1/\e_1,\cdots,x_d/\e_d, t/\delta^2),
	\end{equation*}
	where $(\e_i)_{1\le i\le d} \in (0,1]^d$ and $\delta \in (0,1]$. Then
	\begin{equation}\label{est.heatLip}
		\norm{\nabla u_\e}_{L^\infty(Q_{1/2})} \le C\norm{\nabla u_\e}_{L^2(Q_{1})},
	\end{equation}
	where $C$ depends only on $d,\Lambda,L$ and $\gamma$.
\end{theorem}

\begin{proof}   It suffices to prove  
\begin{align}\label{n-heatLip}
|\nabla u_\e(0,0)| \leq C\Big(\fint_{Q_1} |\nabla u_\va|^2\Big)^{1/2}.
\end{align}
	The proof follows the same idea of Section \ref{sec.reperiodization}, combined with the Lipschitz estimate \eqref{est.Lip.e-d} obtained previously for the case $A_\e(x,t) = A(x,t,x/\e,t/\delta^2)$. Precisely, let $\lambda_i = \e_1/\e_i \ge 1$ for $ 1\leq i\leq d$, and let  $M_\lambda, \Phi_\lambda$ be defined as in \eqref{Mlm} and \eqref{est.PhiPositivity} respectively. Then $v_\e(x,t) = u_\e(\Phi^{-1}_\lambda x,t)$ satisfies
	\begin{equation*}
		\frac{\partial}{\partial t} v_\e-\txt{div} (A^\sharp_\e \nabla v_\e) = 0 \quad \text{in } \Phi_\lambda(B_1)\times(-1,0),
	\end{equation*}
	where $A^\sharp_\e = A^\sharp( x,t,x/\e_1,t/\de^2 )= \Phi_\lambda A(\Phi^{-1}_\lambda x, t, \fl{M_\lm}x/\e_1,  t/\de^2) \Phi_\lambda$. Note that $ A^\sharp (x,t,y,s)$ satisfies the same conditions as $A(x,t,y,s)$. By \eqref{est.Lip.e-d} (or
	\cite[Theorem 1.1]{GN21}, it holds that
	\begin{equation}\label{gnre1}
		\sup_{\delta \vee \e_1 \leq r < 1} \Big(\fint_{Q_r} |\nabla v_\va|^2\Big)^{1/2}
		\le C  \Big(\fint_{Q_1} |\nabla v_\va|^2 \Big)^{1/2}.
	\end{equation}
	As a consequence,  
	\begin{equation} \label{gnre2}
		\sup_{\delta \vee \e_1 \leq r <1} \Big(\fint_{Q_r} |\nabla u_\va|^2\Big)^{1/2}
		\le C  \Big(\fint_{Q_1} |\nabla u_\va|^2 \Big)^{1/2}.
	\end{equation} 
	
	If $\delta \leq \e_d$, by setting $w_\e(x,t)=u_\e(\lm x,\lm^2t) $ with $\lm=\e_1,\e_2\cdots, \e_d,$ sequentially, and finally $\lm=\de$ if $\de<\e_d$, we can perform the rescalling arguments and use \eqref{gnre2} repeatedly as in the proof of Theorem \ref{thm.main2}  to obtain that 
	\begin{equation} \label{n-gnre2}
	\begin{split}
		&\sup_{\e_{j+1} \leq r<\e_{j} }\Big(\fint_{Q_r} |\nabla u_\va|^2\Big)^{1/2}
		\le C  \Big(\fint_{Q_{\e_j}} |\nabla u_\va|^2 \Big)^{1/2}\quad \text{for } j=1,2,\cdots, d-1,\\
		&\sup_{\de \leq r< \e_{d} }\Big(\fint_{Q_r} |\nabla u_\va|^2\Big)^{1/2}
		\le C  \Big(\fint_{Q_{\e_d}} |\nabla u_\va|^2 \Big)^{1/2} \quad\text{ if } \de<\e_d.
		\end{split}
	\end{equation}
	As a result,  
	\begin{equation} \label{gnre3}
		\sup_{\de\wedge \e_d \leq r<1 }\Big(\fint_{Q_r} |\nabla u_\va|^2\Big)^{1/2}
		\le C  \Big(\fint_{Q_1} |\nabla u_\va|^2 \Big)^{1/2}, 
	\end{equation} where $a\wedge b=\min\{a,b\} $.

Likewise, in the case $\de \ge \e_1$, we can perform the rescalling arguments sequentially by setting $w_\e(x,t)= u_\e(\lm x,\lm^2t)$ with $\lm=\de,\e_1,\e_2,\cdots,\e_d$ to derive \eqref{gnre3}. And in the case     $\e_{j_0-1} \ge \delta > \e_{j_0}$ for some $1< j_0 < d$, we can also prove   \eqref{gnre3} by setting $w_\e(x,t)= u_\e(\lm x,\lm^2t)$ with  $\lm=\e_1,\e_2,\cdots, \e_{j_0}, \de, \e_{j_0+1}, \cdots \e_d,$  and performing the rescalling arguments subsequently.

Note that the desired estimate \eqref{n-heatLip} is a direct consequence of \eqref{gnre3} and a standard blow-up argument. We therefore complete the proof.
\end{proof}

\subsection{Stable approximation}
In Subsection \ref{sec.reperiodization}, we have seen that the uniform Lipschitz estimate for the operator $-\txt{div} (A_\e \nabla u_\e)$ holds and is stable under the variable-separation condition \eqref{Var Separation}. In this subsection, we explore the existence and stability of effective approximate problems, and the corresponding error estimate. These stability properties are important in numerical computation for realistic problems. On the other hand, they also provide an alternative approach for Lipschitz estimate or even higher-order regularity in the presence of correctors.

Let us consider the following Dirichlet problem:
 \begin{align} \label{ap-eq-1}
	 -\text{div} (A(x,x_1/\e_1, \cdots,x_d/\e_d) \na u_\va)=F  \quad \text{in } \Omega,  \qquad  u_\va=g\quad\text{on } \pa\Omega.
\end{align}
The next theorem provides the existence of effective approximate problem as well as the corresponding convergence rate, which depends continuously on $(\e_i)_{1\le i\le d}$. In other words,  small changes of $(\e_i)_{1\le i\le d}$ will only lead to a small change on the convergence rate.
\begin{theorem}\label{error-th}
Let $\Omega$ be a bounded $C^{1,1}$ domain. Suppose $A(x,y)$ satisfies conditions  \eqref{ellipticity}, \eqref{periodicity}, and \eqref{ccon}   with $\gamma=1$. Let $u_\e$ be the weak solution to  \eqref{ap-eq-1},
and $u_0$ the solution to
\begin{align}\label{ap-eq-ho}
	 -\txt{div} (\widehat{A^\lm}(x ) \na u _0) =F \quad \text{in } \Omega,
	 \qquad u_0=g  \quad\text{on } \pa \Omega,
\end{align}  where  $\widehat{A^\lm}$ is the effetive matrix given by \eqref{alamhat} below with $\lm=(1, \e_1/\e_2,\cdots, \e_1/\e_d)$. Then
 \begin{align}\label{error-th-1}
 \|  u_\e-u_0\|_{L^2(\Omega)}
 \le
C \{\e_1+\e_2+\cdots+\e_d\}
\| u_0\|_{H^2(\Omega)}
 ,
 \end{align}
 where $C$ depends only on $d$, $\Lambda,$ $ L$ and $\Omega$.
\end{theorem}

We recall that under the scale-separation condition \eqref{cond.separate} the effective approximate solution $u_0$ was found by reiterated homogenization \cite{nsxjfa2020} (including the one-scale case, i.e., $\e_1 = \e_2 = \cdots = \e_d$). However, since condition \eqref{cond.separate} is not assumed and we have general multiscales $(\va_i)_{1\le i\le d}$, the reiterated homogenization method seems not work here. 
In order to find an effective approximate problem (or effective matrix), we introduce a family of parameters $(\lambda_i)_{1\le i\le d}$ to reduce the problem \eqref{ap-eq-1} to a locally periodic problem considered in e.g. \cite{nsxjfa2020}. The payoff is that now, different from \cite{nsxjfa2020}, the periodic cell could be very thin, or by a change of variable, the cell problem may have degenerate coefficients. Fortunately, this could be handled by the idea of reperiodization or the strategy in \cite{GN21}. We point out that the assumption $\gamma = 1$ is not essential and only accounts for the optimality of convergence rate.

\begin{proof}[Proof of Theorem \ref{error-th}.]
For $ \lm=(\lm_1,\cdots,\lm_d)\in (0,+\infty)^d$ , define $M_\lambda$ as in \eqref{Mlm}. We consider a family of $\lm$-dependent operators
\begin{align}\label{Llm}
	\mathcal{L}_\va^\lm = -\text{div} ( A^\lm (x, x/\va_1 ) \na ),
\end{align}
where $A^\lm(x,y) = A(x,M_\lm y)$.
It is obvious that  $A^\lm$ satisfies the ellipticity condition \eqref{ellipticity}, and $A^\lm$ is $(\lm_1^{-1},\cdots,  \lm_{d}^{-1})$-periodic in $(y_1 ,\cdots, y_d )$. 
Let $\chi^\lm (x,y)= (\chi^\lm _{j}(x, y)), 1\leq j\leq d,$ be the correctors given by 
\begin{equation}\label{chilm}
	\left\{
	\aligned
	&- \text{div} ( A^\lm(x,y)\na_y   \chi^\lm _{j} (x,y)  )
	=   \text{div}  ( A^\lm (x,y) e_j )
	~~    \text{ in } \R^{d},\\
	&\chi^\lm _{j}(x, y)  \ \text{is } (\lm_1^{-1},\lm_2^{-1},\cdots,\lm_{d}^{-1}) \text{-periodic in }  (y_1 ,\cdots,y_d ) , \\
	&\int_{0}^{\lm_d^{-1}}\cdots   \int_0^{\lm_1^{-1}} \chi^\lm _{j} (x, y   )\, dy_1  \cdots dy_d  =0,
	\endaligned
	\right.
\end{equation}
where $e_j=(0,...,1,...0)\in \R^d$ with $1$ in the $j$-th position.  For each fixed $\lm $, the homogenized operator of  $\mathcal{L}_\va^\lm$ is given by $ -\text{div} (\widehat{A^\lm } (x) \na )$ with 
\begin{align}\label{alamhat}
	\widehat{A^\lm}(x)=\fint_0^{\lm_d^{-1} }\! \cdots  \fint_0^{ \lm_1^{-1}}A^\lm (x,y) \big(I+  \na_{y}\chi^\lm (x,y )\big )dy_1 \cdots dy_d.
\end{align}
Note that in the above definitions \eqref{chilm} and \eqref{alamhat}, the periodic cell is degenerate if $\lm_i$'s are very large. However, we emphasize that $\widehat{A^\lm}(x)$ is not degenerate and satisfies good properties (such as the ellipticity condition) independent of $\lm_i$'s. To see this, we may use the trick of reperiodization and observe that $A^\lm(x,y)$ is $\Phi_\lm^{-1}$-periodic in $y$, i.e., $A^\lm(x,y) = A^\lm(x,y+\Phi_\lm^{-1}z)$ for any $z\in \Z^d$. Therefore, the periodic cell is a rectangular region with all side lengths in $[\frac12,1]$ and the correctors and homogenized operators are perfectly constructed with non-degenerate estimates.

Let $u_\e^\lm $ and  $u^\lm_0$ be the weak solutions to
\begin{align} \label{ap-eq-va}
	 -\text{div} (A^\lm(x,x/\e_1 ) \na u^\lm_\va)=F  \,\,\,  \text{in } \Omega,  \quad \text{and} \quad u^\lm_\va=g\,\,\, \text{on } \pa\Omega,
\end{align}
 and
\begin{align}\label{ap-eq-ho1}
	 -\text{div} (\widehat{A^\lm}(x ) \na u^\lm_0) =F \quad \text{ in } \Omega,
	\quad \text{and} \quad u^\lm_0=g  \text{ on } \pa \Omega,
\end{align} 
respectively.
By either repeating the aguments in Section 7 of \cite{GN21} (where simialr problem was considered for the parabolic operators) or employing the idea of reperiodization, one can prove that
\begin{align} \label{leap-2re}
		\begin{split}
			\|u^\lm_\va-u_0^\lm\|_{L^2(\Omega)}
			&\leq  C\sum_{i=1}^d \lm_i^{-1}\e_1 \|u_0^\lm\|_{H^2(\Omega)}
		\end{split}
	\end{align}
for some positive constant $C$ depending only on $d,\Lambda$, and $\Omega$. Finally, since $$-\text{div} (A^\lm(x,x/\e_1 ) \na ) = -\text{div} (A (x,x_1/\e_1, \cdots, x_d/\e_d) \na )$$ 
when  $\lm_i=\va_1/\va_{i} , i=1,2, \cdots, d,$
the desired estimate \eqref{error-th-1} follows directly from \eqref{leap-2re}. 
 \end{proof}

The proof above implies that for each $\lm$ the solution $u_\va^\lm$ of problem \eqref{ap-eq-va} could be approximated by an effective solution $u_0^\lm$ of problem \eqref{ap-eq-ho1}.  
This provides us an alternative proof of Theorem \ref{thm.main2} based on the excess decay method \cite{armstrongan2016}, see also \cite{armstrongcpam2016,shenapde2017, ShenBook18}. 
To see this briefly, suppose that $A=A(x, y)$ satisfies \eqref{ellipticity}, \eqref{periodicity} and \eqref{ccon}.  
	Let $u^\lm_\e$ be a weak solution of $ -\text{div}( A^\lm(x,x/\e_1)\na u_\va^\lm)=0$ in $B_{2r}$,
	and set $\de= \max_{1\le i\le d} \lm_i^{-1} \va_1  \le r\le 1$. By performing the same analysis as Section 3 in \cite{GN21} (see \cite{GS20} and \cite{nsxjfa2020} for the original ideas), one can prove that 
	there exists a weak solution $u_0^\lm$ to $ -\text{div} (\widehat{A^\lm} \na u^\lm_0)=0$ in $B_r$ such that
	\begin{equation}\label{5.21}
		\aligned
		& \left(\fint_{B_r} |u^\lm_\e -u^\lm_0 |^2 \right)^{1/2}
		  \le C \left\{ \left(\frac{\de}{r}\right)^\sigma
		+\de^\gamma L  \right\}
		\left(\fint_{B_{2r}} |u^\lm_\e|^2\right)^{1/2}
		 ,
		\endaligned
	\end{equation}
	where  $\sigma>0$  depends only on $d$ and $\Lambda$, and the constant $C$ depends only on $d$, $\Lambda$ and $(L,\gamma)$ in \eqref{ccon}. 
 Since  $u_0^\lm$ is $C^{1,\alpha}$ continuous (note that $\widehat{A^\lm}$ is H\"{o}lder continuous), an application of the excess decay method then gives 
\begin{align}\label{reth41}
\left(\fint_{B_r} |\na u^\lm_\va |^2\right)^{1/2}\leq C  \left(\fint_{B_1}|  \nabla u^\lm_\va |^2\right)^{1/2}  
\end{align}for any $\max_{1\le i\le d} \lm_i^{-1} \va_1  \le r< 1$, where $C$ is independent of $\va$  and $\lm$. 
By taking 
$\lm=(1, \va_1/\va_2,\cdots, \va_1/\va_d)$, and then performing the rescaling argument as we did in Section 4.1, we obtain \eqref{thm.main2-re}.

To end this section, we show the stability of the effective coefficient matrix $\widehat{A^\lm}$ (given by \eqref{alamhat}) with respect to the parameters $\lm = (\lm_i)_{1\le i\le d}$. This tells us that the effective approximate problem  \eqref{ap-eq-ho}  is stable and relatively small changes of $\lambda$ will only lead to a small change of $u_0^\lm$ in the energy sense.
\begin{proposition}
\label{leconAlm}
	For any $
	\lm=(\lm_1, \lm_2,\cdots ,\lm_d), \kappa=(\kappa_1,\kappa_2,\cdots ,\kappa_d) \in (0,+\infty)^d$,  we have 
	\begin{align}\label{leconAlm-re}
		\big\|\widehat{A^{\lm}}-\widehat{A^{\kappa}} \big\|_{\infty} \leq C      \sum_{i=1}^d\big( \kappa_i^{-1} +\lm_i^{-1} \big)  \big|\lm_i-\kappa_i \big|  ,
	\end{align}
	where $C$ depends only on $d$ and $\Lambda$.
\end{proposition}

Let $M_\lm$ be the $d\times d$ diagonal matrix defined as \eqref{Mlm},
  and set  $ \widetilde{\chi}^{\lm }(x,y)=\chi^{\lm}(x, M_\lambda^{-1} y).$
By \eqref{chilm}, $\widetilde{\chi}^\lm _{j} (x,y )$ is $1$-periodic in $y$ and satisfies
\begin{align}\label{chilm1}
	&- \text{div} (M_\lm  A (x,y) M_\lm \na_y  \widetilde{\chi}^\lm_{j} (x,y)  )
	=   \text{div}  ( M_\lm  A  (x,y) e_j )
\quad \text{ in } \mathbb{T}^{d}=[0,1)^d.
\end{align}
The standard energy estimate implies that 
\begin{align}\label{eschilm1}
	\fint_{\mathbb{T}^d}  |M_\lm \na_y \widetilde{\chi}^\lm (x,y) |^2  dy\leq C
\end{align}
for some constant $C$ depending only on $d $ and $\Lambda.$ Now, Proposition \ref{leconAlm} is a consequence of the fact
\begin{align} \label{N-Alamda}
\widehat{A^\lm}(x) = \fint_{\mathbb{T}^d } A  (x,y) \big(I+  M_\lm\na_{y}\widetilde{\chi}^\lm (x,y )\big )dy,\end{align}
and the following lemma.
\begin{lemma}\label{leschi1chi2}
	Let $\widetilde{\chi}^{\lm }$ and $ \widetilde{\chi}^{\kappa }$  be  given by \eqref{chilm1} with respect to the parameter $\lm=(\lm_1, \lm_2, \cdots ,\lm_d),$ $\kappa=(\kappa_1,\kappa_2,\cdots ,\kappa_d) \in (0,+\infty)^d$.   Then
	\begin{align}\label{lechi1chi2-re}
		\begin{split}
			& \int_{\mathbb{T}^d} \big|M_{\lm }\na_y  \widetilde{\chi}^{\lm }(x,y)- M_\kappa \na_y\widetilde{\chi}^{\kappa}(x,y) \big|^2dy
			 \leq   C     \sum_{i=1}^d\big( \kappa_i^{-2} +\lm_i^{-2} \big)  \big|\lm_i-\kappa_i \big|^2 ,
		\end{split}
	\end{align}
	where  $C$ depends on $d $ and $\Lambda$ .
\end{lemma}

\begin{proof}
 By the coerciveness of $A$, for any $1\leq j\leq d$ 
	\begin{align}\label{plechi1chi2-0}
		\begin{split}
			&\Lambda \int_{\mathbb{T}^d}\big|M_\lm \na_y \widetilde{\chi}_j^{\lm }(x,y)-M_\kappa \na_y \widetilde{\chi}_j^{\kappa}(x,y)\big|^2dy\\
			&\leq \int_{\mathbb{T}^d} A \big(M_\lm \na_y \widetilde{\chi}_j^{\lm }(x,y)-M_\kappa \na_y \widetilde{\chi}_j^{\kappa}(x,y)\big)    \cdot \big(M_\lm \na_y \widetilde{\chi}_j^{\lm }(x,y)-M_\kappa \na_y \widetilde{\chi}_j^{\kappa}(x,y) \big)dy.
		\end{split}
	\end{align}
	In view of  \eqref{chilm1},
	\begin{align}\label{eq.chi.l-k}
		\begin{split}
			&-\text{div}\big\{ M_{\lm }A  \big( M_{\lm } \na_y \widetilde{\chi}_j^{\lm }(x,y)-M_{\kappa}\na_y\widetilde{\chi}_j^{\kappa}(x,y)\big)\big\} \\
			&= \text{div} \big\{ \big(M_{\lm }- M_{\kappa}\big) A  M_{\kappa} \na_y \widetilde{\chi}_j^{\kappa}(x,y)  \big\}+\text{div}\big\{ M_{\lm }Ae_j-M_{\kappa}Ae_j\big\}.
		\end{split}
	\end{align}
	Multiplying the above equation with $ \widetilde{\chi}_j^{\lm }(x,y ) $ and  integrating by parts, we arrive at
	\begin{align}\label{plechi1chi2-1}
		\begin{split}
			& \int_{\mathbb{T}^d}  A \big(M_{\lm } \na_y  \widetilde{\chi}_j^{\lm }(x,y) -M_{\kappa} \na\widetilde{\chi}_j^{\kappa}(x,y) \big) \cdot  M_\lm \na_y \widetilde{\chi}_j^{\lm}(x,y)  dy \\
			&=-\int_{\mathbb{T}^d}  A  M_{\kappa}   \na_y  \widetilde{\chi}_j^{\kappa}(x,y)   \cdot \big(M_{\lm }-M_{\kappa}\big)  \na_y \widetilde{\chi}_j^{\lm }(x,y) dy
		- \int_{\mathbb{T}^d}A e_j \cdot \big(M_\lm -M_\kappa \big) \na_y \widetilde{\chi}_j^{\lm }(x,y) dy .
		\end{split}
	\end{align}
	Likewise,  integrating \eqref{eq.chi.l-k} against $\widetilde{\chi}^\kappa_j(x,y)$ leads to
	\begin{align}\label{plechi1chi2-2}
		\begin{split}
			& \int_{\mathbb{T}^d}  A \big(M_\lm \na_y \widetilde{\chi}_j^{\lm}(x,y) -M_\kappa \na_y\widetilde{\chi}_j^{\kappa}(x,y) \big)   \cdot M_\kappa \na_y \widetilde{\chi}_j^{\kappa}(x,y)  dy \\
			&=-\int_{\mathbb{T}^d}  A  M_{\lm }   \na_y  \widetilde{\chi}_j^{\lm }(x,y)   \cdot \big(M_{\lm }-M_\kappa\big)  \na_y \widetilde{\chi}_j^{\kappa}(x,y) dy  
			 -\int_{\mathbb{T}^d}A e_j \cdot \big(M_{\lm } -M_\kappa \big) \na_y \widetilde{\chi}_j^{\kappa}(x,y) dy .
		\end{split}
	\end{align}
	Subtracting \eqref{plechi1chi2-2} from \eqref{plechi1chi2-1}, and using \eqref{plechi1chi2-0}, we get
	\begin{align}\label{plechi1chi2-3}
		\begin{split}
			&   \Lambda \int_{\mathbb{T}^d}\big|M_{\lambda }\na_y \widetilde{\chi}_j^{\lambda }(x,y)-M_\kappa\na_y \widetilde{\chi}_j^{\kappa}(x,y) \big|^2dy\\
			&\leq -\int_{\mathbb{T}^d}   A  M_\kappa   \na_y  \widetilde{\chi}_j^{\kappa} (x,y)   \cdot \big(I_d-(M_{\lambda })^{-1}M_\kappa\big) \big( M_{\lambda } \na_y \widetilde{\chi}_j^{\lambda }(x,y)- M_\kappa \na_y \widetilde{\chi}_j^{\kappa}(x,y)\big)  dy \\
&\quad +\int_{\mathbb{T}^d}   A  \big( M_{\lambda } \na_y \widetilde{\chi}_j^{\lambda }(x,y)- M_\kappa \na_y \widetilde{\chi}_j^{\kappa}(x,y)\big) \cdot \big(I_d-(M_{\lambda })^{-1}M_\kappa\big) M_\kappa   \na_y  \widetilde{\chi}_j^{\kappa} (x,y)  dy \\
			&\quad-\int_{\mathbb{T}^d}   A  M_{\lambda }   \na_y  \widetilde{\chi}_j^{\lambda }(x,y)    \cdot \big( 2I_d-M_{\lambda }(M_\kappa)^{-1}- (M_{\lambda })^{-1}M_\kappa \big) M_\kappa \na_y \widetilde{\chi}_j^{\kappa}(x,y)  dy \\
			&\quad-\int_{\mathbb{T}^d}   A e_j \cdot \big(I_d -(M_{\lambda })^{-1}M_\kappa \big) \big( M_{\lambda }\na_y\widetilde{\chi}_j^{\lambda }(x,y) - M_\kappa\na_y \widetilde{\chi}_j^{\kappa}(x,y) \big)  dy\\
           &\quad-\int_{\mathbb{T}^d}  A e_j \cdot \big(2I_d -M_{\lambda }(M_\kappa)^{-1}-(M_{\lambda })^{-1}M_\kappa \big) M_\kappa\na_y \widetilde{\chi}_j^{\kappa}(x,y)  dy\\
           &\doteq E_1+E_2+\cdots+E_5.
		\end{split}
	\end{align}
	Note that
 	\begin{align*}
		&\big| I_d- (M_\lm)^{-1} M_\kappa\big|\leq C  \max_{1\leq i \leq d}  \big\{\kappa_i^{-1}  |\lm_i-\kappa_i |\big\},\\
		&\big| I_d- (M_\kappa)^{-1} M_\lm\big|\leq C \max_{1\leq i \leq d}  \big\{\lm_i^{-1}  |\lm_i-\kappa_i |\big\}.
	\end{align*}
By \eqref{eschilm1} and the Cauchy-Schwarz inequality,  we deduce that 
\begin{align*}
 E_1+E_2+E_4 \leq C  \max_{1\leq i \leq d}  \big\{ (\kappa_i^{-2}+\lm_i^{-2}) |\lm_i-\kappa_i |^2\big\}  + \frac{\Lambda}{4} \int_{\mathbb{T}^d}\big|M_{\lambda }\na_y \widetilde{\chi}_j^{\lambda }(x,y)-M_\kappa\na_y \widetilde{\chi}_j^{\kappa}(x,y) \big|^2dy.
 \end{align*}
In view of \eqref{eschilm1} and the fact that
 \begin{align*}
		\big| 2I_d- (M_\lm)^{-1} M_\kappa- M_{\lambda }(M_\kappa)^{-1}\big|\leq C  \max_{1\leq i \leq d}  \big\{(\kappa_i \lm_i)^{-1}   |\lm_i-\kappa_i |^2\big\},
	\end{align*}
we can bound $E_3$ and $ E_5$ as following
\begin{align*}
 E_3+E_5  \leq C  \max_{1\leq i \leq d}  \big\{ (\kappa_i\lm_i)^{-1} |\lm_i-\kappa_i |^2\big\}.  
 \end{align*}
By substituting the estimates of $E_1$-$E_5$ into \eqref{plechi1chi2-3},	one derives \eqref{lechi1chi2-re}  immediately.
\end{proof}

\begin{remark}
From \eqref{chilm1} and \eqref{N-Alamda}, we observe that the effective matrix $\widehat{A^\lm}$ is scaling invariant, namely, $\widehat{A^\lm} = \widehat{A^{t\lm}}$ for any $t >0$. Therefore, in the statement of Theorem \ref{error-th}, we can choose $\lambda = (1/\e_1, 1/\e_2,\cdots, 1/\e_d)$, which is symmetric in $(\e_i)_{1\le i\le d}$. Also,
\eqref{leconAlm-re} can be strengthened to 
\begin{align*}
		\big\|\widehat{A^{\lm}}-\widehat{A^{\kappa}} \big\|_{\infty} \leq C \inf_{t>0}  \   \sum_{i=1}^d\big( \kappa_i^{-1} +(t \lm_i)^{-1} \big)  \big|t\lm_i-\kappa_i \big|.
	\end{align*}
	Particularly, the right-hand side of the above inequality is zero if $\kappa = t_0 \lm$ for some $t_0 >0$.
\end{remark}

\section{Reiterated homogenization of quasi-periodic operators}\label{sec.quasipH}
This section is devoted to the proof of Theorem \ref{thm.quasiper.HC},  i.e., the reiterated homogenization of $-\text{div}(A(x,x/\e_1,\cdots,x/\e_n)\na )$ under the assumption that $A(x,y_1,\cdots,y_n)$ is quasi-periodic in $y_i\in \R^d, 1\leq i \leq n$, which means that
there exist a matrix function  $B(x,w_1,w_2,\cdots, w_n)$, $1$-periodic in $w_i \in \R^{m_i}$ with $d\le m_i\in \mathbb{N}$ for $ 1\leq i \leq n$, i.e.,
\begin{equation}\label{b-peri}
	B(x,w_1+z_1, \cdots, w_n+z_n) = B(x,w_1, \cdots,w_n),    ~~~\txt{for any}~ w_i\in \R^{m_i}  , z_i\in \Z^{m_i},  
\end{equation}
 and  matrices 
 $ M_i\in \mathcal{M}^{m_i\times d}$ (the set of $m_i\times d$ constant real matrices) satisfying
\begin{align} \label{nondeg}
	M^T_i z\neq 0  \quad\text{ for any} \quad  z\in \Z^{m_i} \setminus \{0 \},\quad 1\leq i \leq n,
\end{align} such that 
$A(x,y_1,\cdots,y_n) = B(x,M_1 y_1, \cdots, M_n y_n).$  
Through out this section, we shall always assume that $(\e_i)_{1\leq i\leq n}$ satisfies the scale-separation condition \eqref{cond.separate}. Our proof follows closely the ideas of \cite{ngestseng,allaire1992,allaire1996} and \cite{guenean2010,guenean2018}.

\subsection{Multiscale convergence}

Let $\Omega$ be a bounded domain in $\R^d$, and let $\mathbb{T}^m=[0,1)^m$ denote the $m$-dimensional torus for any $m\in \mathbb{N}$. 
Consider the Hilbert space $L^2 (\Omega \times \mathbb{T}^{m_1}\times\cdots \times \mathbb{T}^{m_{n}} )$, equipped with the inner product
\begin{equation*}
    \Ag{f,g} = \int_{\Omega} \int_{\mathbb{T}^{m_1}}\cdots \int_{\mathbb{T}^{m_n}} f(x,y_1,\cdots, y_n) \cdot \overline{g(x,y_1,\cdots,y_n)} \ dx dy_1\cdots dy_n.
\end{equation*}
For any $f\in L^2 (\Omega \times \mathbb{T}^{m_1}\times\cdots \times \mathbb{T}^{m_{n}} )$, we can write $f$ in terms of Fourier series
\begin{equation}
    f(x,y_1,\cdots,y_n) = \sum_{\substack{k_i\in \mathbb{Z}^{m_i},  1\le i\le n}} \lambda(k_1,\cdots,k_n;x) \exp\Big( 2\pi \mathbf{i} \sum_{1\le j\le n} k_j\cdot y_j\Big),
\end{equation}
and the norm of $f$ is given by
\begin{equation}\label{def.L2norm}
    \| f\|_2^2:= \Ag{f,f} = \sum_{\substack{k_i\in \mathbb{Z}^{m_i},  1\le i\le n}} \int_{\Omega} |\lambda(k_1,\cdots, k_n;x)|^2 dx < \infty.
\end{equation}
 \begin{definition}
	Let $M_i\in \mathcal{M}^{m_i\times d}, 1\leq i\leq n,$ satisfy \eqref{nondeg}.
	A family of functions $u_\e (x)\in L^2(\Omega)$ is said to multiscale cut-and-project (associated with $M_i, 1\leq i\leq n$) converge to a function $u_0(x,y_1,\cdots ,y_n) \in L^2(\Omega\times \mathbb{T}^{m_1}\times\cdots \times \mathbb{T}^{m_n})$ if
	\begin{align*}
		\begin{split}
			&\lim_{\e\rightarrow 0} \int_\Omega u_\e(x) \psi(x,\frac{M_1x}{\e_1},\cdots ,\frac{M_nx}{\e_n})dx\\
			&=\int_\Omega\int_{\mathbb{T}^{m_1}}\cdots \int_{\mathbb{T}^{m_n}} u_0(x,y_1,\cdots ,y_n) \psi(x,y_1,\cdots ,y_n)dx dy_1 \cdots  dy_n,
		\end{split}
	\end{align*}
	for any $\psi(x,y_1,\cdots ,y_n)\in L^2(\Omega; C_{per}(\mathbb{T}^{m_1}\times\cdots  \times \mathbb{T}^{m_n})).$
\end{definition}

Note that the above multiscale cut-and-project convergence is stronger than the usual weak convergence in $L^2(\Omega)$, but weaker than the usual strong convergence in $L^2(\Omega)$. The following theorem shows that any bounded sequence in $L^2(\Omega)$ has a multiscale cut-and-project convergent subsequence.

\begin{theorem}\label{quasi-th1}
Let $\{u_\e\}$ be a uniformly bounded sequence in $L^2(\Omega)$. Then up to a subsequence, $u_\e$ multiscale cut-and-project (associated with $M_i, 1\leq  i\leq n$, satisfying \eqref{nondeg}) converges to a limit $u_0(x,y_1,\cdots,y_n)\in L^2(\Omega\times \mathbb{T}^{m_1}\times\cdots  \times \mathbb{T}^{m_n})$.
\end{theorem}
The proof of the above theorem relies on the next two lemmas. For convenience, from now on, for any $\phi(x,y_1,\cdots, y_n) \in L^2(\Omega; C_{per}(\mathbb{T}^{m_1}\times\cdots  \times \mathbb{T}^{m_n}))$, we define
\begin{equation}
    [\phi]_\e(x) := \phi (x, \frac{M_1x}{\e_1},\cdots ,\frac{M_nx}{\e_n}),
\end{equation}
and
\begin{equation}
    [\phi](x) =  \int_{\mathbb{T}^{m_1}}\cdots \int_{\mathbb{T}^{m_n}}  \vp(x, y_1,\cdots ,y_n)dy_1 \cdots d y_n.
\end{equation}
These definitions also apply if $\phi$ is independent of $x$ or some of $y_i$'s.
\begin{lemma}\label{quasi-le1}
	  Let  $\varphi\in C_{per}(\mathbb{T}^{m_1} \times\cdots \times \mathbb{T}^{m_n})$.   Then $[\varphi]_\e(x) = \varphi (\frac{M_1x}{\e_1},\cdots ,\frac{M_nx}{\e_n})$ converges weak-$\star$ in $L^\infty(\Omega)$ to the constant $[\varphi]$.
\end{lemma}
\begin{proof}
	We first consider the trigonometric polynomial functions
	$$[\vp]_\e (x)=\sum_{\substack{ k_i\in \Z^{m_i}, i=1,\cdots ,n\\ |k_i|\leq K_0}} C(k_1,\cdots ,k_n) \exp{ \Big( 2\pi   \mathbf{i} \sum_{j=1}^n k_j\cdot \frac{M_j x}{\e_j} } \Big),$$
	where $K_0 > 0$ and $ C(k_1,\cdots ,k_n)$ denote the Fourier coefficients depending on $ k_1,\cdots ,k_n.$  Since $(\e_i)_{1\le i\le n}$ satisfies the condition \eqref{cond.separate}, by \eqref{nondeg} we know that 
	\begin{equation*}
	    \exp{ \Big( 2\pi \mathbf{i} \sum_{j=1}^n k_j\cdot \frac{M_j x}{\e_j} } \Big) = \exp{ \Big(2 \pi \mathbf{i} \sum_{j=1}^n M_j^T k_j\cdot \frac{x}{\e_j} } \Big)
	\end{equation*}
	converges weakly to the mean value zero \cite{guenean2010}, if there exists some $k_i\in \Z^{m_i}\setminus \{0\}$. As a result,  we have that for any $\Psi(x)\in C(\Omega)\cap  L^1(\Omega)$,
	\begin{align}\label{prquasi-le1-01}
		\lim_{\e\rightarrow 0}\int_\Omega [\vp]_\e(x) \Psi(x) dx = C(0,\cdots,0) \int_\Omega \Psi(x) dx= [\varphi] \int_\Omega \Psi(x) dx.
	\end{align}
	This proves the lemma for trigonometric polynomials.
	Finally, since the trigonometric polynomials are dense in $C_{per}(\mathbb{T}^{m_1} \times\cdots \times \mathbb{T}^{m_n})$, the general case follows by a standard density argument.
	\end{proof}

\begin{lemma}\label{quasi-le2}
	  Let $\vp \in L^2(\Omega; C_{per}(\mathbb{T}^{m_1} \times\cdots \times \mathbb{T}^{m_n} ))$. Then
	\begin{align*}
		\lim_{\e \rightarrow 0}\int_\Omega \big|[\vp]_\e(x) \big|^2  dx=    \int_\Omega\int_{\mathbb{T}^{m_1}}\cdots \int_{\mathbb{T}^{m_n}}   \big|\vp(x,y_1,\cdots ,y_n)\big|^2 dx dy_1 \cdots  dy_n .
	\end{align*}
\end{lemma}
\begin{proof}
	We first consider the case $\vp(x,y_1,\cdots ,y_n)=\al(x)\be(y_1,\cdots ,y_n)$, where $\al(x)  \in L^\infty(\Omega)$ and $ \be(y_1,\cdots ,y_n) \in C_{per}(\mathbb{T}^{m_1}\times \cdots  \times \mathbb{T}^{m_n}).$
	Since $ \be^2  \in C_{per}(\mathbb{T}^{m_1} \times  \cdot \cdot\cdot \times \mathbb{T}^{m_n} )$,  Lemma \ref{quasi-le1} implies that
	$ [\be^2]_\e $ converges weak-$\star$ in $L^\infty(\Omega)$ to $[\be ^2]$.
	As a result,
	\begin{align*}
		\lim_{\e\to 0} \int_\Omega \big|[\vp]_\e(x) \big|^2 dx & =  [\beta^2] \int_\Omega \al^2(x) dx\\
		& =\int_\Omega\int_{\mathbb{T}^{m_1}}\cdots \int_{\mathbb{T}^{m_n}} \Big|\vp (x, y_1,\cdots ,y_n)\Big|^2 dx dy_1\cdots  dy_n.
	\end{align*}
    where we have used the Fubini's theorem in the second identity.
	
One can then extend the above result to step functions 
	\begin{align*}
		\vp_k =\sum_{i=1}^k t_i \chi_{A_i}(x)\psi_i(y_1,\cdots ,y_n),
	\end{align*}
	where $A_i$ are measurable sets in $\Omega$, $\chi_{A_i}$ are the characteristic functions of $A_i$ and $\psi_i(y_1,\cdots y_n)\in C_{per}(\mathbb{T}^{m_1}\times \cdots   \times\mathbb{T}^{m_n}).$
	Finally, by the density of step functions in $L^2(\Omega; C_{per}(\mathbb{T}^{m_1}\times \cdots \times   \mathbb{T}^{m_n})),$ one derives the desired result.
\end{proof}

\begin{proof}[\textbf{Proof of Theorem \ref{quasi-th1} }]
	Let $\{u_\e\}$ be a bounded sequence in $L^2(\Omega)$. Consider the linear form
	\begin{align*}
		\langle \mathrm{L}_\e ,\vp\rangle=\int_\Omega [\vp]_\e(x) u_\e(x)dx, \qquad \vp\in L^2(\Omega; C_{per}(\mathbb{T}^{m_1}\times \cdots \times  \mathbb{T}^{m_n})).
	\end{align*}
	The Cauchy inequality
	 implies that
	$\mathrm{L}_\e$  is a bounded sequence of continuous linear maps on  $L^2(\Omega; C_{per}(\mathbb{T}^{m_1}  \times\cdot \cdot\cdot \times \mathbb{T}^{m_n} )).$ Therefore there exists a subsequence of $(\e_i)_{1\leq i \leq n}$, still denoted by $(\e_i)_{1\leq i \leq n}$, and a linear map $\mathrm{L}$ in the dual of $L^2(\Omega; C_{per}(\mathbb{T}^{m_1}  \times\cdot \cdot\cdot \times \mathbb{T}^{m_n} ))$ such that
	\begin{align*}
		\langle \mathrm{ L} ,\vp\rangle  = \lim_{ \va\rightarrow 0} \langle \mathrm{L}_{\e},\vp\rangle,
	\end{align*} for any $\vp \in  L^2(\Omega; C_{per}(\mathbb{T}^{m_1}  \times\cdot \cdot\cdot \times \mathbb{T}^{m_n} )).$
Thanks to Lemma \ref{quasi-le2},
	\begin{align*}
		|\langle  \mathrm{L} ,\vp\rangle| \leq  C \lim_{ \va\rightarrow 0} \big\|  [\vp]_\e  \big\|_{L^2(\Omega)}  \leq C \|\vp\|_{L^2(\Omega\times \mathbb{T}^{m_1}\times\cdots \times \mathbb{T}^{m_n})}.
	\end{align*}
	Thus $\mathrm{L}$ can be extended to a bounded linear functional on $L^2(\Omega\times \mathbb{T}^{m_1}\times\cdots  \times\mathbb{T}^{m_n})$. By the Riesz representation theorem, there exists a  $u_0\in L^2(\Omega\times \mathbb{T}^{m_1}\times\cdots  \times \mathbb{T}^{m_n})$ such that
	\begin{align*}
		\lim_{\va\rightarrow 0}\int_\Omega [\vp]_\e(x) u_\e(x)dx = \int_\Omega\int_{\mathbb{T}^{m_1}}\cdots \int_{\mathbb{T}^{m_n}} u_0(x,y_1,\cdots ,y_n) \vp(x,y_1,\cdots ,y_n)dx dy_1\cdots  dy_n.
	\end{align*}
	The proof is complete.
\end{proof}

\subsection{Multiscale convergence of gradients}
In this subsection, we will deal with the multiscale cut-and-project convergence of $\nabla u_\e$. 

For each $1\leq \ell \leq n$, let $M_\ell\in \mathcal{M}^{m_\ell\times d}$  satisfy \eqref{nondeg}. Define
\begin{align*}
	\mathcal{L}_{\ell}&=\Big\{   w \in L^2 (\Omega \times \mathbb{T}^{m_1}\times\cdots \times \mathbb{T}^{m_{n}} )^d \big| \\
	& \qquad w = \sum_{\substack{k_i\in \mathbb{Z}^{m_i},  1\le i\le\ell  \\ k_\ell \neq 0}} \lm(k_1,\cdots ,k_\ell;x) M^T_\ell k_\ell \exp\Big( 2\pi \mathbf{i} \sum_{1\le j\le \ell} k_j\cdot y_j\Big)    \Big\},
\end{align*}
where $\lambda(k_1,\cdots,k_n;x)$ are scalar functions.

Define
\begin{equation*}
    \mathcal{H}_n = \big\{  w\in  L^2 (\Omega \times \mathbb{T}^{m_1}\times\cdots \times \mathbb{T}^{m_{n}} )^d \big|  (M_n^T\na_{y_n})\cdot w=0 \big\},
\end{equation*}
and for $1\le \ell \le n-1$,
\begin{align*}
    \mathcal{H}_\ell &= \big\{  w\in  L^2 (\Omega \times \mathbb{T}^{m_1}\times\cdots \times \mathbb{T}^{m_{n}} )^d \big|  
		(M^T_\ell \na_{y_\ell})\cdot\int_{\mathbb{T}^{m_{\ell+1}}     }\cdots \int_{\mathbb{T}^{m_n}     }  w=0 \big\}.
\end{align*}
Let
\begin{equation}
    \mathcal{H} = \bigcap_{1\le \ell \le n} \mathcal{H}_\ell \quad \txt{and} \quad \mathcal{L} = \sum_{1\le \ell \le n} \mathcal{L}_\ell.
\end{equation}
It is clear that $\mathcal{H}_\ell$ and $\mathcal{H}$ are all closed subspaces of $L^2 (\Omega \times \mathbb{T}^{m_1}\times\cdots \times \mathbb{T}^{m_{\ell}} )^d$.

\begin{proposition} \label{prop.LH}
$\mathcal{L} = \mathcal{L}_1 \oplus \mathcal{L}_2 \oplus\cdots \oplus \mathcal{L}_n$ and $L^2 (\Omega \times \mathbb{T}^{m_1}\times\cdots \times \mathbb{T}^{m_{n}} )^d = \mathcal{L} \oplus \mathcal{H}$.
\end{proposition}

\begin{proof}
We first show that $\mathcal{L}_\ell$ is closed in $L^2 (\Omega \times \mathbb{T}^{m_1}\times\cdots \times \mathbb{T}^{m_{n}} )^d$. Let $\{ w^{(m)}: m\ge 1 \}$ be a Cauchy sequence in $\mathcal{L}_\ell$ that converges to some $w^* \in L^2 (\Omega \times \mathbb{T}^{m_1}\times\cdots \times \mathbb{T}^{m_{n}} )^d$. We would like to show $w^* \in \mathcal{L}_\ell$.

Write
\begin{equation*}
    w^{(m)} = \sum_{\substack{k_i\in \mathbb{Z}^{m_i},  1\le i\le\ell  \\ k_\ell \neq 0}} \lm^{(m)}(k_1,\cdots ,k_\ell;x) M^T_\ell k_\ell \exp\Big( 2\pi \mathbf{i} \sum_{1\le j\le \ell} k_j\cdot y_j\Big).
\end{equation*}
By Parseval's identity, for each $(k_1,\cdots, k_\ell)$ with $k_\ell \neq 0$, $\{ \lm^{(m)}(k_1,\cdots, k_\ell; x)M_\ell^T k_\ell : m\ge 1 \}$ is a Cauchy sequence in $L^2(\Omega)^d$. Since $|M^T_\ell k_\ell|>0$ is fixed, $\{ \lm^{(m)}(k_1,\cdots, k_\ell; x) : m\ge 1 \}$ is a Cauchy sequence in $L^2(\Omega)$. Let $\lm^*(k_1,\cdots, k_\ell;\cdot)$ be the limit of $\lm^{(m)}(k_1,\cdots, k_\ell;\cdot)$ as $m\to \infty$. Then, $\{ \lm^{(m)}(k_1,\cdots, k_\ell; x)M_\ell^T k_\ell : m\ge 1 \}$ converges to $ \lm^{*}(k_1,\cdots ,k_\ell;x) M^T_\ell k_\ell $ and by the uniqueness of the limit, we must have
\begin{equation*}
    w^{*} = \sum_{\substack{k_i\in \mathbb{Z}^{m_i},  1\le i\le\ell  \\ k_\ell \neq 0}} \lm^{*}(k_1,\cdots ,k_\ell;x) M^T_\ell k_\ell \exp\Big( 2\pi \mathbf{i} \sum_{1\le j\le \ell} k_j\cdot y_j\Big).
\end{equation*}
This implies that $w^* \in \mathcal{L}_\ell$ and thus $\mathcal{L}_\ell$ is closed. 
Note that $\mathcal{L}_\ell, 1\le \ell \le n,$ are mutually orthogonal. We conclude that  $\mathcal{L} = \mathcal{L}_1\oplus \cdots \oplus \mathcal{L}_n$, and therefore $\mathcal{L}$ is also closed.

Next, we show that $\mathcal{L}_\ell^\perp= \mathcal{H}_\ell$. We first prove that $\mathcal{H}_\ell\subset\mathcal{L}_\ell^\perp $.   Let $v \in \mathcal{H}_\ell$. Consider $k_\ell \neq 0$ and
\begin{equation}\label{eq.Ll.w}
    w = \lambda(k_1,\cdots,k_\ell; x) M^T_{\ell} k_\ell \exp\Big( 2\pi \mathbf{i} \sum_{1\le j\le \ell} k_j\cdot y_j\Big).
\end{equation}
Since $M^T_{\ell} k_\ell \neq 0$, we can write $w = M^T_\ell \nabla_{y_\ell} \widetilde{w}$, where
\begin{equation}
    \widetilde{w} = \frac{1}{2\pi \mathbf{i}} \lambda(k_1,\cdots,k_\ell; x) \exp\Big( 2\pi \mathbf{i} \sum_{1\le j\le \ell} k_j\cdot y_j\Big).
\end{equation}
It follows from the definition of $\mathcal{H}_\ell$ and the fact that $\widetilde{w}$ depends only on $y_1,\cdots, y_\ell$ that
\begin{equation}
    \Ag{w,v} = \Ag{M^T_\ell \nabla_{y_\ell} \widetilde{w}, v  } = -\int_{\Omega} \int_{\mathbb{T}^{m_1}}\cdots \int_{\mathbb{T}^{m_n}} \widetilde{w} \big( M^T_\ell \nabla_{y_\ell}\cdot v \big) = 0.
\end{equation}
Thus, $w\perp \mathcal{H}_\ell$. By considering the linear combinations of $w$'s in the form \eqref{eq.Ll.w} and using a density argument, we see that $\mathcal{H}_\ell \subset\mathcal{L}_\ell^\perp $. 

Now, we prove that $\mathcal{L}_\ell^\perp  \subset \mathcal{H}_\ell $. Consider an arbitrary $v\in \mathcal{L}_\ell^\perp $. Then $\Ag{w,v} = 0$ for any $w$ in the form of \eqref{eq.Ll.w}. Write $v$ as
\begin{equation}
    v = \sum_{\substack{k_i\in \mathbb{Z}^{m_i},  1\le i\le n}} \alpha(k_1,\cdots,k_n;x) \exp\Big( 2\pi \mathbf{i} \sum_{1\le j\le n} k_j\cdot y_j\Big).
\end{equation}
Then
\begin{equation}
    \Ag{w ,v} = \int_{\Omega} \lambda(k_1,\cdots,k_\ell; x) M^T_{\ell} k_\ell \cdot \alpha(k_1,\cdots,k_\ell,0\cdots,0;x) dx = 0.
\end{equation}
Since the above identity holds for any function $\lambda(k_1,\cdots,k_n;\cdot) \in L^2(\Omega)$, by duality, we  have 
\begin{equation}\label{eq.MTlkla}
    M^T_{\ell} k_\ell \cdot \alpha(k_1,\cdots,k_\ell,0\cdots,0;x) = 0
\end{equation}
By taking different $w$'s in \eqref{eq.Ll.w} with $k_\ell \neq 0$, we have \eqref{eq.MTlkla} for every $(k_1,\cdots, k_\ell)$ with $k_\ell\neq 0$. This implies
\begin{equation*}
    (M^T_\ell \na_{y_\ell})\cdot\int_{\mathbb{T}^{m_{\ell+1}}     }\cdots \int_{\mathbb{T}^{m_n}     }  v=0,
\end{equation*}
in the sense of distribution (with an obvious modification if $\ell = n$). Therefore, $v\in \mathcal{H}_\ell$ and thus $\mathcal{L}_\ell^\perp \subset \mathcal{H}_\ell$. Combined with $\mathcal{H}_\ell \subset \mathcal{L}_\ell^\perp$ proved previously, this leads to $\mathcal{L}_\ell^\perp = \mathcal{H}_\ell$.

Finally, since $\sum_{1\le \ell \le n}  \mathcal{H}_\ell^\perp=\sum_{1\le \ell \le n}  \mathcal{L}_\ell=\mathcal{L} $ is closed, a classical result in functional analysis implies that 
\begin{equation}
  \mathcal{H}^\perp= \Big(\bigcap_{1\le \ell\le n} \mathcal{H}_\ell \Big)^\perp= \sum_{1\le \ell \le n}  \mathcal{H}_\ell^\perp    = \mathcal{L},
\end{equation} 
and therefore, 
 $L^2 (\Omega \times \mathbb{T}^{m_1}\times\cdots \times \mathbb{T}^{m_{n}} )^d = \mathcal{L} \oplus \mathcal{H}$.  This completes the proof.
 \end{proof}

The following is the main theorem of this subsection.
\begin{theorem}\label{quasi-th3}
	Let $\{u_\e (x)\}$ be a  bounded sequence in $H^1_0(\Omega)$, and let $M_i\in \mathcal{M}^{m_i\times d}, 1\leq i \leq n,$ satisfy \eqref{nondeg}. Then there are $u_0\in L^2(\Omega)$ and $U_i \in \mathcal{L}_i$ such that up to subsequences,
	\begin{align}
		&u_\e \longrightarrow u_0, \quad\,\,\,\,
		 \nabla u_\e  \longrightarrow \nabla_x u_0+ \sum_{i=1}^n U_i,
	\end{align}
	in the sense of multiscale cut-and-project convergence.
\end{theorem}

It is important to point out that, in view of the definition of $\mathcal{L}_i$, the function $U_i$ can formally be written as $M^T_i \nabla_{y_i} u_i(x,y_1,\cdots, y_i)$. However, $u_i $ may not be in $L^2 (\Omega \times \mathbb{T}^{m_1}\times\cdots \times \mathbb{T}^{m_{\ell}} )$, because $M_i^T k_i $ (though not equal to 0) could be very singular as $k_i$ increases to infinity. This is the key difference between periodic and quasi-periodic cases.

Let $\Xi$ be the subspace of  $ \mathcal{D} (\Omega; C^\infty_{per}(\mathbb{T}^{m_1}\times\cdots \times \mathbb{T}^{m_n}) )^d \cap \mathcal{H}$ composed of partial Fourier series with respect to $(y_1,\cdots ,y_n)$, namely,
\begin{align*}
	\Xi&= \Big\{v\in  \mathcal{D} (\Omega; C^\infty_{per}(\mathbb{T}^{m_1}\times\cdots \times \mathbb{T}^{m_n}) ) ^d \Big| v=\!\!\!\!\sum_{\substack{k_i\in \mathbb{Z}^{m_i},\\  i=1,\cdots ,n,\\ |k_1|+\cdots+|k_n|\leq N}}\!\!\!\!\Lm(k_1,\cdots ,k_n;x)  \exp\Big( 2\pi \mathbf{i} \sum_{1\le j\le n} k_j\cdot y_j\Big) \\
	&  \text{ for some } N\in \mathbb{N} , \text{ and }   \Lm(k_1,\cdots ,k_\ell,0,\cdots,0; x) \cdot M^T_\ell k_\ell=0  \,\,\forall x\in \Omega,  \text{ if  } k_\ell \neq 0,  1\leq \ell \leq n \Big\}.
\end{align*}
Clearly, $\Xi$ is dense in $\mathcal{H}.$ Similarly, for $\ell\in {1,2,\cdots ,n}$, let $F_\ell$ be the subspace of $\mathcal{D}(\Omega; C^\infty_{per}(\mathbb{T}^{m_1}\times\cdots \times \mathbb{T}^{m_\ell}   ))$ composed of partial Fourier series with respect to $(y_1,\cdots ,y_\ell)$ with zero mean in $y_\ell$.

The proof of Lemma \ref{quasi-th3} relies essentially on the following lemma.
\begin{lemma}\label{h-1bound}
	Let $M_i\in \mathcal{M}^{m_i\times d}, 1\leq i \leq n,$ satisfy \eqref{nondeg}.
	  Then for any $\varphi \in F_\ell$,
\begin{align}\label{th-h-1bound-re}
		\frac{1}{\va_\ell} [\varphi]_\va(x) =\frac{1}{\va_\ell} \varphi\Big(x, \frac{M_1x}{\va_1}, \cdots ,\frac{M_\ell x}{\va_\ell} \Big) \quad \text{is uniformly (in $\e$) bounded in } H^{-1}(\Omega).
		\end{align}
\end{lemma}

\begin{corollary} \label{co1}
	Let $\varphi$ belong  to the subspace of $\mathcal{D}(\Omega; C^\infty_{per}(\mathbb{T}^{m_1}\times\cdots \times \mathbb{T}^{m_n}))$ composed of partial Fourier series, such that
	\begin{align*}
		\int_{\mathbb{T}^{m_\ell}   }\cdots   \int_{\mathbb{T}^{m_n}} \varphi dy_\ell\cdots dy_n =0 \quad \text{ for some } 1 \le \ell  \le n.
	\end{align*}
	Then $\frac{1}{\va_\ell} [\varphi]_\va$ is uniformly (in $\e$) bounded in $H^{-1}(\Omega)$.
\end{corollary}

\begin{proof} The proof is the same as Corollary 3.4 in \cite{allaire1996}. The key idea is to rewrite $ \frac{1}{\va_\ell} [\varphi ]_\va$ as
	\begin{align*}
		\frac{1}{\va_\ell} [\varphi ]_\va = \sum_{j=\ell}^n\frac{\va_j}{\va_\ell} \frac{1}{\va_j} [\varphi_j]_\va,
	\end{align*}
	where $\varphi_n=\varphi-\int_{\mathbb{T}^{m_n}}\varphi,$ and \begin{align*}
	 			\varphi_j=\int_{\mathbb{T}^{m_{j+1}} }\cdots \int_{\mathbb{T}^{m_n}} \varphi\, dy_{j+1}\cdots dy_n -\int_{\mathbb{T}^{m_{j}}  }\cdots \int_{\mathbb{T}^{m_n}} \varphi\, dy_{j}\cdots dy_n
	\end{align*} for $\ell \leq j\leq n-1$, satisfy the assumption of Lemma \ref{h-1bound}. Therefore, $\e_j^{-1} [\varphi_j]_\e$ is uniformly bounded in $H^{-1}(\Omega)$ and so is $\e_\ell^{-1} [\varphi]_\e$.
\end{proof}
\begin{proof}[\textbf{Proof of Lemma \ref{h-1bound}}]
    Since $F_\ell$ contains the partial Fourier series, to show \eqref{th-h-1bound-re}, it is sufficient to consider a single term in the form of
	\begin{align} \label{def.varphi}
	\varphi=  \Lm(k_1,\cdots ,k_\ell;x)  \exp\Big( 2\pi \mathbf{i} \sum_{1\le j\le \ell} k_j\cdot y_j\Big)
	\end{align} 
	for some $k_\ell \neq 0$ and $\Lm(k_1,\cdots,k_\ell;x) \in L^\infty(\Omega)$. In the following, we will use the infinity norm $\| f \|_\infty: = \| f \|_{L^\infty(\Omega\times \mathbb{T}^{m_1}\times \cdots \times \mathbb{T}^{m_\ell})}$, because of the property $\| [f]_\e \|_\infty \le \| f\|_{\infty}.$
	
	Define 
	$$ S\varphi:=   \frac{ M^T_\ell k_\ell}{2\pi \mathbf{i}|M^T_\ell k_\ell|^2}\Lm(k_1,\cdots ,k_\ell;x)  \exp\Big( 2\pi \mathbf{i} \sum_{1\le j\le n} k_j\cdot y_j\Big) = \frac{ M^T_\ell k_\ell}{2\pi \mathbf{i}|M^T_\ell k_\ell|^2} \varphi.$$
	It is easy to see that
	\begin{align}\label{h-1bound-1}
		(M^T_\ell \na_{y_\ell})\cdot S\varphi=\varphi. 
	\end{align}
	For any given value of $(a_1,\cdots ,a_{\ell-1})\in [0,1]^{\ell-1}$, let $T$ be a linear operator from $F_\ell$ to $F_\ell$ given by
	\begin{align}\label{def.T}
		T\varphi:= \sum_{j=1}^{\ell-1} a_j (M^T_j\na_{y_j}) \cdot S  \varphi = \theta \varphi,
	\end{align}
	where $\theta$ is a constant scalar given by
	\begin{equation*}
	    \theta = \sum_{j=1}^{\ell-1} a_j (M^T_jk_j) \cdot \frac{ M^T_\ell k_\ell}{|M^T_\ell k_\ell|^2}.
	\end{equation*}
	Note that $|\theta| \le C_0$ for some $C_0$ depending only on $M_j$ and $k_j$ with $1\le j\le \ell$.

    By the previous definition of $S$ and $T$, it is obvious that $T^m \varphi = \theta^m \varphi$ and $ST^m \varphi = \frac{ M^T_\ell k_\ell}{2\pi \mathbf{i}|M^T_\ell k_\ell|^2} \theta^m \varphi$ for any integer $m \ge 1$. Thus,
	\begin{align}\label{h-1bound-8}
		\| T^m \varphi \|_\infty \le C_0^m \| \varphi \|_{\infty}  \quad \txt{and} \quad \| S T^m \varphi \|_\infty \le C_1 C_0^m \| \varphi \|_{\infty},
	\end{align}
	where $C_1 = |\frac{ M^T_\ell k_\ell}{2\pi \mathbf{i}|M^T_\ell k_\ell|^2}|$ depends only on $M_\ell$ and $k_\ell$.
	
	Now, we are ready to prove \eqref{th-h-1bound-re} for the particular $\varphi$ given by \eqref{def.varphi}. Since $(\e_i)_{1\le i\le \ell}$ satisfies the scale-separation condition \eqref{cond.separate}, there exists $\e^*>0$ such that if $\e_\ell < \e^*$ (Recall that $(\e_i)_{1\le i\le \ell}$ is in decreasing order)
	\begin{align}\label{h-1bound-9}
		 r_\e:= C_0 \frac{\va_\ell}{\va_{\ell-1}} \leq \frac{1}{2}.
	\end{align}
	Note that
	\begin{align*}
		\frac{1}{\va_\ell}[\varphi]_\va=\text{div}[S\varphi]_\va-[\text{div}_xS\varphi]_\va-\frac{\va_\ell}{\va_{\ell-1}}\frac{1}{\va_\ell}[T_\va\varphi]_\va,
	\end{align*}
	where
	$
		T_\va  \varphi=\sum_{j=1}^{\ell-1}\frac{\va_{\ell-1}}{\va_j} (M_j^T\na_{y_j})\cdot S\varphi.
	$ Observe that $T_\e$ takes the same form of \eqref{def.T} with $a_j = \e_{\ell-1}/\e_j \in [0,1]$.
	Iterating the above equality $m$ times yields
	\begin{align}\label{h-1bound-10}
		\begin{split}
			\frac{1}{\va_\ell}[\varphi]_\va&=\sum_{p=0}^{m-1}(-1)^p \Big(\frac{\va_\ell}{\va_{\ell-1}}\Big)^p \big( \text{div}[ST_\va^p \varphi]_\va- [\text{div}_x (ST_\va^p\varphi)]_\va\big)\\
			&\quad+(-1)^m \Big(\frac{\va_\ell}{\va_{\ell-1}}\Big)^m\frac{1}{\va_\ell}[T_\va^m\varphi]_\va.
		\end{split}
	\end{align}
	Note that for a smooth vector-valued function $\Psi(x,y_1,\cdots,y_\ell)$ periodic in $(y_1,\cdots,y_\ell)$,
	\begin{align}\label{est.Psi.e}
		\|[\Psi]_\va\|_{L^2(\Omega)}  +\|\text{div}[\Psi]_\va\|_{H^{-1}(\Omega)}\leq C(\Omega)  \| \Psi\|_\infty.
	\end{align}
	In view of \eqref{h-1bound-10} and \eqref{est.Psi.e}, we have
	\begin{align*}
		\begin{split}
			& \frac{1}{\va_\ell} \|[\varphi]_\va\|_{H^{-1}(\Omega)}\\
			&\leq C(\Omega) \Big\{\frac{1}{\va_\ell}\Big(\frac{\va_\ell}{\va_{\ell-1}}\Big)^m \|T_\va^m \varphi\|_\infty +  \sum_{p=0}^{m-1} \Big(\frac{\va_\ell}{\va_{\ell-1}}\Big)^p \big(\|ST_\va^p \varphi\|_\infty + \|\text{div}_x (ST_\va^p\varphi)\|_\infty \big)\Big\}\\
			&\leq C(\Omega) \Big\{\frac{1}{\va_\ell}\Big(\frac{\va_\ell}{\va_{\ell-1}}\Big)^m C_0^m\| \varphi\|_\infty +  \sum_{p=0}^{m-1} \Big(\frac{\va_\ell}{\va_{\ell-1}}\Big)^p C_1 C_0^{p}\big(\| \varphi\|_\infty+ \|\na_x\varphi\|_\infty \big)\Big\},
		\end{split}
	\end{align*}
	where we have used \eqref{h-1bound-8} and the fact $\| [f]_\e \|_\infty \le \| f\|_{\infty}$.  In view of \eqref{h-1bound-9}, we derive that
	\begin{align*}
		\begin{split}
			\frac{1}{\va_\ell} \|[\varphi]_\va\|_{H^{-1}(\Omega)}
			&\leq C(\Omega) \Big\{\frac{1}{\va_\ell}(r_\va)^m \|\varphi\|_\infty   + C_1 (\|\varphi\|_\infty+ \|\na_x\varphi\|_\infty)
			\sum_{p=0}^{m-1} (r_\va)^p\Big\}\\
			&\leq   C(\Omega) \Big\{\frac{1}{\va_\ell}(r_\va)^m \|\varphi\|_\infty    +\frac{C_1}{1-r_\e}(\|\varphi\|_\infty+ \|\na_x\varphi\|_\infty) \Big\}.
		\end{split}
	\end{align*}\
	Since $r_\va\leq1/2$, by letting $m$ tend to infinity we get
	\begin{equation*}
	    \frac{1}{\va_\ell} \|[\varphi]_\va\|_{H^{-1}(\Omega)} \le 2C(\Omega)C_1 (\|\varphi\|_\infty+ \|\na_x\varphi\|_\infty),
	\end{equation*}
	for any $\e_\ell < \e^*$. Finally, for $\e_\ell \ge \e^*$, a trivial bound gives
	\begin{equation*}
	    \frac{1}{\va_\ell} \|[\varphi]_\va\|_{H^{-1}(\Omega)} \le \frac{1}{\e^*} C(\Omega) \|\varphi \|_\infty.
	\end{equation*}
	Combining the last two estimates for different ranges of $\e_\ell$, we obtain \eqref{th-h-1bound-re}.
\end{proof}

\begin{proof}[\textbf{Proof of Theorem \ref{quasi-th3}}]
	Since  $\{u_\e\}$ is a bounded sequence in $H_0^1(\Omega)$, by Theorem \ref{quasi-th1}, there exist two functions $u_0(x,y_1,\cdots ,y_n)\in L^2(\Omega\times \mathbb{T}^{m_1}\times\cdots \times \mathbb{T}^{m_n})$ and $  \Theta_0(x,y_1,\cdots ,y_n) \in L^2(\Omega\times \mathbb{T}^{m_1}\times\cdots \times \mathbb{T}^{m_n})^d$ such that $u_\e$ and $\na u_\e$  multiscale cut-and-project converge (up to subsequences) to $u_0$ and $\Theta_0$, respectively.
	
	We first show that $u_0$ is independent of $y_1,\cdots, y_n$. By the convergence of $\na u_\e$,  for $ \Psi \in  \mathcal{D}(\Omega; C^\infty_{per}(\mathbb{T}^{m_1}\times\cdots \times \mathbb{T}^{m_n} ))^d$,
	\begin{align} \label{prquasi-th3-01}
		\lim_{\e \rightarrow 0} \int_{\Omega} \na u_\e  (x) \cdot [\Psi]_\e(x) dx=\int_\Omega\int_{\mathbb{T}^{m_1}}\cdots \int_{\mathbb{T}^{m_n}} \Theta_0 \cdot\Psi dx dy_1\cdots dy_n .
	\end{align}
	Observe that
	\begin{equation}\label{eq.div.Psie}
	    \txt{div}( [\Psi]_\e) = [\txt{div}_x \Psi]_\e + \sum_{i=1}^n \frac{1}{\e_i} [M^T_i \nabla_{y_i}\cdot \Psi]_\e.
	\end{equation}
	Consequently, by the integration by parts,
\begin{align}\label{def.Ie}
		\begin{split}
			I_\e &:= \int_{\Omega} \na u_\e  (x) \cdot [\Psi]_\e(x) dx\\
			&= -\int_{\Omega}  u_\e  (x) [\txt{div}_x \Psi]_\e(x)dx  - \int_{\Omega}  u_\e (x) \sum_{i=1}^n\frac{1}{\e_i} [M^T_i \nabla_{y_i}\cdot \Psi]_\e(x) dx.
		\end{split}
	\end{align}	
	
	In view of \eqref{prquasi-th3-01}, we see that
	\begin{align}\label{prquasi-th3-02}
	\lim_{\e\rightarrow 0} \e_n I_\e=	\lim_{\e\rightarrow 0}\e_n \int_\Omega\int_{\mathbb{T}^{m_1}}\cdots \int_{\mathbb{T}^{m_n}} \Theta_0\cdot \Psi dx dy_1\cdots dy_n =0.
	\end{align}
	On the other hand, since $(\e_i)_{1\le i\le n}$ satisfies the scale-separation condition \eqref{cond.separate}, we have
	\begin{align}
		&\lim_{\e\rightarrow 0}\e_n \int_{\Omega}  u_\e (x) [\text{div}_x\Psi]_\e(x) dx = 0  ,\label{eq.limit.Ie}\\
	&\lim_{\e\rightarrow 0}	\e_n \int_{\Omega}  u_\e (x) \sum_{i=1}^{n-1}\frac{1}{\e_i} [M^T_i \nabla_{y_i}\cdot \Psi]_\e(x) dx =0  .\label{eq.limit.Ie1}
	\end{align}
	The equality \eqref{def.Ie}, combined with \eqref{prquasi-th3-02}, \eqref{eq.limit.Ie} and \eqref{eq.limit.Ie1}, implies that
	\begin{align*}
		  \lim_{\e\rightarrow 0}\int_{\Omega}  u_\e (x)   [M^T_n \nabla_{y_n}\cdot \Psi]_\e(x) dx=0.
			\end{align*}
	By the multiscale cut-and-project convergence of $u_\e$  and integration by parts,
	\begin{align*}
		\begin{split}
		 0&=  \int_\Omega\int_{\mathbb{T}^{m_1}}\cdots \int_{\mathbb{T}^{m_n}}u_0(x,y_1,\cdots, y_n)     M_n^T \na_{y_n}\cdot \Psi (x,y_1,\cdots ,y_n )dx dy_1 \cdots  dy_n\\
			&=-\int_\Omega\int_{\mathbb{T}^{m_1}}\cdots \int_{\mathbb{T}^{m_n}} M_n^T \na_{y_n}u_0(x,y_1,\cdots ,y_n)\cdot    \Psi (x,y_1,\cdots ,y_n )dx dy_1 \cdots  dy_n,
		\end{split}
	\end{align*}
	for any $\Psi \in  \mathcal{D}(\Omega; C^\infty_{per}(\mathbb{T}^{m_1}\times\cdots \times \mathbb{T}^{m_n} ))^d$. Thus, in the sense of distribution,
	\begin{align}\label{prquasi-th3-05}
		M_n^T \na_{y_n}u_0(x,y_1,\cdots ,y_n) =0.
	\end{align}
	Since $u_0\in L^2(\Omega\times \mathbb{T}^{m_1} \times\cdots \times \mathbb{T}^{m_n} )$, we can write
	\begin{align*}
		u_0(x,y_1,\cdots,y_n)=\sum_{ k_i\in \mathbb{Z}^{m_i}} C(x,k_1,\cdots ,k_n) \exp\Big( 2\pi \mathbf{i} \sum_{1\le j\le n} k_j\cdot y_j\Big),
	\end{align*}
	which together with \eqref{prquasi-th3-05}, implies that
	\begin{align*}
		\sum_{ k_i\in \mathbb{Z}^{m_i}}2\pi \mathbf{i} M_n^Tk_n C(x,k_1,\cdots ,k_n) \exp\Big( 2\pi \mathbf{i} \sum_{1\le j\le n} k_j\cdot y_j\Big)=0.
	\end{align*}
	Thus, we have $
		2\pi \mathbf{i}  M_n^Tk_n C(x,k_1,\cdots ,k_n)=0.	$
	Since $ M_n^Tk_n\neq 0$ for any $k_n\in \Z^{m_n}\setminus \{0\} $, $C(x,k_1,\cdots,k_n)=0$ for any $k_n\neq 0$ . It follows that  $u_0$ is independent of $y_n$. Performing similar argument repeatedly, we can show that $u_0$ is independent of $y_{n-1},\cdots ,y_1$. Hence, $u_0=u_0(x)$ depends only on $x$.
	
	To identify $ \Theta_0$, we take $ \Psi \in  \Xi$. By \eqref{eq.div.Psie}, integration by parts and the fact $(M^T_n \na_{y_n}) \cdot \Psi=0$, we have
	\begin{align}\label{prquasi-th3-07}
		\int_\Omega \na u_\va \cdot [\Psi]_\va  dx =- \int_\Omega u_\va [\text{div}_x \Psi]_\va dx -\sum_{\ell=1}^{n-1}\frac{1}{\va_\ell}\int_\Omega u_\va [ (M^T_\ell\na_{y_\ell}) \cdot \Psi]_\e.
	\end{align} By the definition of $ \Xi$, $(M^T_\ell\na_{y_\ell}) \cdot \Psi$ satisfies the assumption of Corollary \ref{co1} with $\ell$ replaced by $\ell+1$. As a result,
	$
	\frac{1}{\va_{\ell+1}} [(M^T_\ell \na_{y_\ell}) \cdot \Psi]_\e
	$ is bounded in  $H^{-1}(\Omega)$, and $ \frac{1}{\va_{\ell}} [(M^T_\ell \na_{y_\ell}) \cdot \Psi]_\e \rightarrow 0$ strongly in  $H^{-1}(\Omega)$ for all $1\leq \ell\leq n-1$ as $\va\rightarrow 0.$ Therefore,  letting $\e  \rightarrow 0$ in \eqref{prquasi-th3-07}, we get
	\begin{align*}
		\begin{split}
			&\int_\Omega\int_{\mathbb{T}^{m_1}}\cdots \int_{\mathbb{T}^{m_n}} \Theta_0(x,y_1,\cdots ,y_n) \cdot \Psi(x,y_1,\cdots ,y_n)dy_1 \cdots  dy_n dx\\
			&=- \int_{\Omega}  \int_{\mathbb{T}^{m_1}}\! \cdots \!\int_{\mathbb{T}^{m_n}} u_0(x) (\text{div}_x\Psi) (x,y_1,\cdots ,y_n ) dy_1 \cdots  dy_n  dx\\
			& = \int_{\Omega}  \int_{\mathbb{T}^{m_1}}\! \cdots \!\int_{\mathbb{T}^{m_n}} \nabla u_0(x) \cdot \Psi (x,y_1,\cdots ,y_n ) dy_1 \cdots  dy_n  dx,
		\end{split}
	\end{align*}
	which implies that $u_0\in H^1(\Omega)$, and for $\Psi \in \Xi$,
	\begin{align*}
			&\int_\Omega\int_{\mathbb{T}^{m_1}}\cdots \int_{\mathbb{T}^{m_n}} (\Theta_0(x,y_1,\cdots ,y_n)-\na u_0) \Psi(x,y_1,\cdots ,y_n) dy_1 \cdots  dy_n dx=0.
	\end{align*}       
	Now, since $\Xi $ is dense in $\mathcal{H}$, the above equation implies $ \Theta_0 -\na u_0 \in \mathcal{H}^\perp.$
	As a consequence of Proposition \ref{prop.LH},
	\begin{align*}
		\na u_\va(x,y_1,\cdots ,y_n) \longrightarrow \Theta_0 =  \na u_0 +\sum_{i=1}^n U_i
	\end{align*} in the sense of multiscale cut-and-project convergence, where $U_i \in \mathcal{L}_i$. This ends the proof of the theorem.
\end{proof}

\subsection{Homogenization theorem}
In this subsection, we prove the reiterated homogenization theorem for multiscale quasi-periodic operators.
\begin{theorem}\label{quasi-th4}
	Let $M_i$ satisfy \eqref{nondeg}, and $B=B(x,y_1,\cdots,y_n) \in L^\infty(\Omega; C_{per}(\mathbb{T}^{m_1}\times\cdots  \times \mathbb{T}^{m_n}))^{d\times d}$ satisfy \eqref{ellipticity}-\eqref{smoothness}.
	Let $u_\va\in H_0^1(\Omega)$ be the solution to $-\txt{div}([B]_\e \nabla u_\e)=f$ in $\Omega$ with $f\in H^{-1}(\Omega)$. Then $u_\e$ converges weakly in $ H_0^1(\Omega)$ to a function $u_0$, and $\na u_\va$ converges in the sense of multiscale cut-and-project convergence to $\na u_0(x) +\sum_{\ell=1}^n U_\ell,$ where $U_\ell \in \mathcal{L}_\ell, 1\leq \ell\leq n,$ are given by the system
	\begin{equation}\label{quasith4-re1}
		\begin{cases}
			-\txt{div}_{y_n}\big(M_n B (\na u_0+ \sum_{\ell=1}^n U_\ell )\big)=0,\\
			-\txt{div}_{y_j}\Big( \int_{\mathbb{T}^{m_{j+1}} } \cdots \int_{\mathbb{T}^{m_n}   }  M_j B \big(\na u_0+ \sum_{\ell=1}^n U_\ell \big) dy_{j+1}\cdots dy_{n}\Big)  =0, \quad 1\leq j\leq n-1, \\
			-\txt{div}_x\Big( \int_{\mathbb{T}^{m_1}   } \cdots \int_{\mathbb{T}^{m_n}   }   B \big(\na u_0+ \sum_{\ell=1}^n U_\ell \big) dy_1\cdots dy_n\Big)  =f(x).
		\end{cases}
	\end{equation}
	Moreover, $u_0$ is the unique solution of $-\txt{div} (B_0(x)\na u_0)=f$ in $ \Omega$ with $B_0(x)$ satisfying \eqref{ellipticity} and \eqref{smoothness}, defined by an inductive formula.
	Precisely, for $0\le k\le n-1$,
	\begin{align}\label{def.Bk}
		B_k = \int_{\mathbb{T}^{m_{k+1}}   } B_{k+1}(I_d + \X_{k+1} ) dy_{k+1},
	\end{align} where $\X_{k+1} = (\X_{k+1, j}) \in \mathcal{L}_{k+1}^{d}$ is the unique solution to
	\begin{align}\label{def.Xk}
		\int_{\mathbb{T}^{m_{k+1}}   } B_{k+1} \big(e_j+ \X_{k+1,j} \big)\cdot M^T_{k+1}\na_{y_{k+1}} \phi\, dy_{k+1}=0, \qquad 1\le j\le d,
	\end{align}
	for any $\phi\in C_{per}^\infty(\T^{m_{k+1}})$.
\end{theorem}
\begin{proof}
	First of all, since $\{u_\e\}$ is bounded in $H^1_0(\Omega)$, there exists $u_0\in H_0^1(\Omega)$ such that, up to a subsequence, $u_\e$ converges to $u_0$ weakly in $H_0^1(\Omega)$ and strongly in $L^2(\Omega)$ (this implies that $u_\e$ converges to $u_0$ in the sense of multiscale cut-and-project convergence). Thanks to Theorem \ref{quasi-th3}, there exist functions $U_i \in \mathcal{L}_i$ ($1\le i\le n$) such that, up to subsequences,
	$ \na u_\va \longrightarrow  \na u_0+\sum_{\ell=1}^n U_\ell
	$ in the sense of multiscale cut-and-project convergence.  Let $\varphi\in \mathcal{D}(\Omega)$ and $\varphi_\ell \in \mathcal{D}\big(\Omega; C^\infty_{per}(\mathbb{T}^{m_1} \times\cdots  \times \mathbb{T}^{m_\ell}   )\big) $ for $1\leq \ell\leq n$. Taking $\varphi+\sum_{\ell=1}^n \va_\ell [\varphi_\ell]_\va$, where $[\varphi_\ell]_\va=\varphi_\ell\big(x, \frac{M_1x}{\e_1},\cdots, \frac{ M_\ell x} {\e_\ell}\big)  \in H^1_0(\Omega)$, as a test function in the variational formulation of $ -\txt{div}([B]_\e \nabla u_\e)=f$, we have
	\begin{align}\label{eq.ue.variational}
		\begin{split}
			\int_\Omega [B]_\va \na u_\va \cdot \na \big(\varphi+\sum_{\ell=1}^n \va_\ell [\varphi_\ell]_\va\big)dx
			 &=\int_\Omega [B]_\va \na u_\va \cdot  \big(\na_x\varphi+ \sum_{\ell=1}^n \sum_{j=1}^\ell \va_\ell (\va_j)^{-1} [M^T_j\na_{y_j}\varphi_\ell]_\va\big)dx\\
			&=\int_\Omega f \big(\varphi+\sum_{\ell=1}^n \va_\ell [\varphi_\ell]_\va\big) dx.
		\end{split}
	\end{align}
	By the scale-separation condition \eqref{cond.separate}, we have $\lim_{\e\to 0} \e_\ell/\e_j =0$ for all $1\le j\le \ell-1$.
	In view of \eqref{eq.ue.variational} and Theorem \ref{quasi-th1} ($[B]_\e$ should be combined with the test functions),  we derive that
	\begin{align}\label{pquasith4-02}
		\int_\Omega \!\int_{\mathbb{T}^{m_1}} \!\!\cdots \!\! \int_{\mathbb{T}^{m_n}}  B(\na u_0 +  \sum_{\ell=1}^n U_\ell  \big)  \big(\na \varphi+\sum_{\ell=1}^n    M^T_\ell\na_{y_\ell}\varphi_\ell \big)dxdy_1 \cdots dy_n=\int_\Omega  f \varphi dx,
	\end{align}
	for all $\varphi\in \mathcal{D}(\Omega)$ and $\varphi_\ell \in \mathcal{D}\big(\Omega; C^\infty_{per}(\mathbb{T}^{m_1} \times\cdots  \times \mathbb{T}^{m_\ell}   )\big).$ Note that the above variational equation also shows that the flux $[B]_\e \nabla u_\e$ multiscale cut-and-project converges to $B(\nabla u_0 + \sum_{\ell=1}^n U_\ell)$.

	Since the functions in the form of $\sum_{\ell=1}^n    M^T_\ell\na_{y_\ell}\varphi_\ell$ are dense in $\mathcal{L}$, the bilinear form in \eqref{pquasith4-02} is bounded and coercive in the Hilbert space
	$$V = H^1_0(\Omega) \times \mathcal{L}, \quad \txt{equipped with norm } \|\nabla u_0\|_{L^2(\Omega)}+ \sum_{\ell=1}^n \| U_\ell \|_{2}.$$
By the Lax-Milgram theorem and the fact $\mathcal{L} = \mathcal{L}_1\oplus \cdots \oplus \mathcal{L}_n$,  the solution $(u_0, U_1, \cdots, U_\ell)$ of \eqref{pquasith4-02} is unique. Note that the above argument is independent of the subsequences selected at the begaining of the proof. Therefore, the whole sequences $u_\va$ and $\na u_\va$ converge, respectively, to $u_0$ and $ \na u_0(x)+\sum_{\ell=1}^n U_\ell$.  

Now, we show \eqref{quasith4-re1}. Let $\varphi = 0$ and $\varphi_{\ell}=0$ for all $1\le \ell \le n-1$. So the only nonzero test function is $\varphi_n$ in the variational form \eqref{pquasith4-02}. This is equivalent to the first equation of \eqref{quasith4-re1}. The remaining equations in \eqref{quasith4-re1} can be obtained similarly by choosing certain nonzero test functions.

	Finally, to isolate the scalar equation for $u_0$ from the system \eqref{quasith4-re1}, we will eliminate $U_\ell$'s through a process of reiterated homogenization. To eliminate the function $U_n$, we construct a $d\times d$ matrix $\X_n = (\X_{n,j}) \in \mathcal{L}_n^d $ such that
	\begin{align}\label{prcoroan-2}
		U_n(x,y_1,\cdots ,y_n)= \X_n(x,y_1,\cdots, y_n) \big(\na u_0+ \sum_{\ell=1}^{n-1} U_\ell \big).
	\end{align}
	Substituting this into the first equation of \eqref{quasith4-re1}, we see that $\X_{n} = (\X_{n,j})_{1\le j\le d}$ satisfies
	\begin{equation}\label{eq.Xnj}
		-\text{div}_{y_n} \big( M_n B (e_j+    \X_{n,j} \big)=0, \quad \txt{for each } 1\le j\le d,
	\end{equation}
	or, equivalently, the variational equation \eqref{def.Xk} with $k = n-1$.
	To see the existence of the solution $\X_{n,j} \in \mathcal{L}_n$, we consider the following regularized equation
	\begin{equation}\label{eq.regularized}
		-\text{div}_{y_n} \big( M_n B M_n^T \nabla_{y_n} \chi_{n,j}^\rho \big) - \rho^2 \Delta_{y_n} \chi_{n,j}^\rho= \text{div}_{y_n} \big( M_n B e_j).
	\end{equation}
	Then, for each $\rho>0$, \eqref{eq.regularized} has a unique solution $\chi_{n,j}^\rho(x,y_1,\cdots, y_{n-1},\cdot)$ in $H^1(\mathbb{T}^{m_n})/\R$. Moreover, the energy estimate implies that $M_n^T \nabla_{y_n} \chi_{n,j}^\rho$ is uniformly bounded in $\mathcal{L}_n$ (independent of $\rho$). Thus, by letting $\rho \to 0$, $M_n^T \nabla_{y_n} \chi_{n,j}^\rho$ converges weakly to some $\X_{n,j} \in \mathcal{L}_n$, which satisfies \eqref{eq.Xnj}. The uniqueness of $\X_{n,j}$ can be derived directly from \eqref{eq.Xnj} or \eqref{def.Xk} and the fact that $\{M_n^T \nabla_{y_n} \phi\, |\, \phi\in C_{per}^\infty(\T^{m_{n}})\} $ is dense in $\mathcal{L}_n$.
	
	Now substituting \eqref{prcoroan-2} into the remaining equations in \eqref{quasith4-re1}, we obtain a new system for $u_0$ and $U_\ell$ with $1\le \ell \le n-1$.
	\begin{equation}\label{quasith4-re2}
		\begin{cases}
			-\text{div}_{y_{n-1}}\big(M_{n-1} B_{n-1} (\na  u_0+ \sum_{\ell=1}^{n-1} U_\ell)\big)=0\\
			-\text{div}_{y_j}\Big( \int_{\mathbb{T}^{m_{j+1}} } \cdots \int_{\mathbb{T}^{m_{n-1}}}  M_j B_{n-1} \big(\na  u_0+ \sum_{\ell=1}^{n-1} U_{\ell} \big) dy_{j+1}\cdots d y_{n-1}\Big)=0, \, 1\leq j\leq n-2 \\
			-\text{div}_x\Big( \int_{\mathbb{T}^{m_1}   } \cdots \int_{\mathbb{T}^{m_{n-1}}}   B_{n-1} \big(\na  u_0+ \sum_{\ell=1}^{n-1} U_\ell \big) dy_1\cdots dy_{n-1}\Big) =f(x),
		\end{cases}
	\end{equation}
	where $B_{n-1}=B_{n-1}(x,y_1,\cdots ,y_{n-1})$ is given by
	\begin{align*}
		B_{n-1}=\int_{\mathbb{T}^{m_n}} B \big( I_d +  \X_n \big) dy_n.
	\end{align*}
	It is obvious that $B_{n-1}$ is bounded and satisfies \eqref{periodicity} and \eqref{smoothness}. Moreover since
	\begin{align*}
		(B_{n-1})_{ij}\xi_i\xi_j&=\int_{\mathbb{T}^{m_n}}B  \big(e_j +  \X_{n,j} \big)\xi_j \cdot \big(e_i +  \X_{n,i} \big)\xi_i  dy_n  \nonumber\\
		&\geq \int_{\mathbb{T}^{m_n}} 	\Lambda  |\xi + \X_n \xi|^2 dy_n
		 \geq 	\Lambda  |\xi|^2,\nonumber
	\end{align*}
	$B_{n-1}$ also satisfies the coercivity in \eqref{ellipticity}. Hence, the system \eqref{quasith4-re2} is analogous to \eqref{quasith4-re1} with one less equation and unknown function. Iterating the process above, we can construct $\X_k$ and $B_k$ satisfying \eqref{def.Xk} and \eqref{def.Bk}, and end up with the equation $-\txt{div}(B_0 \nabla u_0) = f$. This completes the proof.
\end{proof}

\begin{proof}[Proof of Theorem \ref{thm.quasiper.HC}] This is a corollary of Theorem \ref{quasi-th4}. In fact, by the definition of quasi-periodic coefficients, $A(x,x/\e_1,\cdots,x/\e_n) = [B]_\e = B(x,M_1x/\e_1,\cdots,M_nx/\e_n)$. Thus, Theorem \ref{quasi-th4} implies that $A(x,x/\e_1,\cdots,x/\e_n)$ $H$-converges to $B_0 = B_0(x)$, which satisfies \eqref{ellipticity} and \eqref{smoothness}.
\end{proof}
\bibliographystyle{amsplain}

\bibliography{NZ2021}

\providecommand{\bysame}{\leavevmode\hbox to3em{\hrulefill}\thinspace}
\providecommand{\MR}{\relax\ifhmode\unskip\space\fi MR }
\providecommand{\MRhref}[2]{%
  \href{http://www.ams.org/mathscinet-getitem?mr=#1}{#2}
}
\providecommand{\href}[2]{#2}
\begin{thebibliography}{10}

\bibitem{allaire1992}
G.~Allaire, \emph{Homogenization and two-scale convergence}, SIAM J. Math.
  Anal. \textbf{23} (1992), no.~6, 1482--1518. \MR{1185639}

\bibitem{allaire1996}
G.~Allaire and M.~Briane, \emph{Multiscale convergence and reiterated
  homogenisation}, Proc. Roy. Soc. Edinburgh Sect. A \textbf{126} (1996),
  no.~2, 297--342. \MR{1386865}

\bibitem{agk16}
N.S. Armstrong, A.~Gloria, and T.~Kuusi, \emph{Bounded correctors in almost
  periodic homogenization}, Arch. Ration. Mech. Anal. \textbf{222} (2016),
  no.~1, 393--426. \MR{3519974}

\bibitem{armstrongcpam2016}
S.~N. Armstrong and Z.~Shen, \emph{Lipschitz estimates in almost-periodic
  homogenization}, Comm. Pure Appl. Math. \textbf{69} (2016), no.~10,
  1882--1923. \MR{3541853}

\bibitem{armstrongan2016}
S.~N. Armstrong and C.~K. Smart, \emph{Quantitative stochastic homogenization
  of convex integral functionals}, Ann. Sci. \'Ec. Norm. Sup\'er. (4)
  \textbf{49} (2016), no.~2, 423--481. \MR{3481355}

\bibitem{al87}
M.~Avellaneda and F.~Lin, \emph{Compactness methods in the theory of
  homogenization}, Comm. Pure Appl. Math. \textbf{40} (1987), no.~6, 803--847.
  \MR{910954}

\bibitem{al89}
\bysame, \emph{Compactness methods in the theory of homogenization. {II}.
  {E}quations in nondivergence form}, Comm. Pure Appl. Math. \textbf{42}
  (1989), no.~2, 139--172. \MR{978702}

\bibitem{lions1978}
A.~Bensoussan, J.-L. Lions, and G.~Papanicolaou, \emph{Asymptotic analysis for
  periodic structures}, AMS Chelsea Publishing, Providence, RI, 2011, Corrected
  reprint of the 1978 original [MR0503330]. \MR{2839402}

\bibitem{bbmm05}
A.~Bondarenko, G.~Bouchitt\'{e}, L.~Mascarenhas, and Rajesh Mahadevan,
  \emph{Rate of convergence for correctors in almost periodic homogenization},
  Discrete Contin. Dyn. Syst. \textbf{13} (2005), no.~2, 503--514. \MR{2152402}

\bibitem{guenean2010}
G~Bouchitt\'{e}, S.~Guenneau, and F.~Zolla, \emph{Homogenization of dielectric
  photonic quasi crystals}, Multiscale Model. Simul. \textbf{8} (2010), no.~5,
  1862--1881. \MR{2769084}

\bibitem{briane2019}
M.~Briane and G.~A. Francfort, \emph{A two-dimensional labile aether through
  homogenization}, Comm. Math. Phys. \textbf{367} (2019), no.~2, 599--628.
  \MR{3936127}

\bibitem{guenuean2019}
E.~Cherkaev, S.~Guenneau, H.~Hutridurga, and N.~Wellander, \emph{Quasiperiodic
  composites: Multiscale reiterated homogenization}, IEEE proceedings of
  Metamaterials, Rome, Italy, (2019), X086--X08.

\bibitem{chipot1986}
M.~Chipot, D.~Kinderlehrer, and G.~Vergara-Caffarelli, \emph{Smoothness of
  linear laminates}, Arch. Rational Mech. Anal. \textbf{96} (1986), no.~1,
  81--96. \MR{853976}

\bibitem{dongjde2019}
R.~Dong, D.~Li, and L.~Wang, \emph{Directional homogenization of elliptic
  equations in non-divergence form}, J. Differential Equations \textbf{268}
  (2020), no.~11, 6611--6645. \MR{4075552}

\bibitem{dongdcds2018}
R.~Dong and L.~Li, D.and~Wang, \emph{Regularity of elliptic systems in
  divergence form with directional homogenization}, Discrete Contin. Dyn. Syst.
  \textbf{38} (2018), no.~1, 75--90. \MR{3708152}

\bibitem{GN21}
J.~Geng and W.~Niu, \emph{Homogenization of locally periodic parabolic
  operators with non-self-similar scales}, arXiv:2103.01418 (2021).

\bibitem{GS20}
J.~Geng and Z.~Shen, \emph{Homogenization of parabolic equations with
  non-self-similar scales}, Arch. Ration. Mech. Anal. \textbf{236} (2020),
  no.~1, 145--188. \MR{4072212}

\bibitem{gloria2019}
A.~Gloria and M.~Ruf, \emph{Loss of strong ellipticity through homogenization
  in 2{D} linear elasticity: a phase diagram}, Arch. Ration. Mech. Anal.
  \textbf{231} (2019), no.~2, 845--886. \MR{3900815}

\bibitem{milton1986}
G.~W. Milton, \emph{Modelling the properties of composites by laminates},
  Homogenization and effective moduli of materials and media ({M}inneapolis,
  {M}inn., 1984/1985), IMA Vol. Math. Appl., vol.~1, Springer, New York, 1986,
  pp.~150--174. \MR{859415}

\bibitem{milton1999}
G.~W. Milton and V.~Nesi, \emph{Optimal {$G$}-closure bounds via stability
  under lamination}, Arch. Ration. Mech. Anal. \textbf{150} (1999), no.~3,
  191--207. \MR{1738117}

\bibitem{ngestseng}
G.~Nguetseng, \emph{A general convergence result for a functional related to
  the theory of homogenization}, SIAM J. Math. Anal. \textbf{20} (1989), no.~3,
  608--623. \MR{990867}

\bibitem{nsxjfa2020}
W.~Niu, Z.~Shen, and Y.~Xu, \emph{Quantitative estimates in reiterated
  homogenization}, J. Funct. Anal. \textbf{279} (2020), no.~11, 108759, 39.
  \MR{4152234}

\bibitem{Oleinik1988}
O.~A. Ole\u{\i}nik, A.~S. Shamaev, and G.~A. Yosifian, \emph{On the
  homogenization of stratified structures}, Analyse math\'{e}matique et
  applications, Gauthier-Villars, Montrouge, 1988, pp.~401--419. \MR{956970}

\bibitem{past07}
S.~E. Pastukhova and R.~N. Tikhomirov, \emph{Operator estimates in reiterated
  and locally periodic homogenization}, Dokl. Akad. Nauk \textbf{415} (2007),
  no.~3, 304--309. \MR{2458607}

\bibitem{shenapde2015}
Z.~Shen, \emph{Convergence rates and {H}\"{o}lder estimates in almost-periodic
  homogenization of elliptic systems}, Anal. PDE \textbf{8} (2015), no.~7,
  1565--1601. \MR{3399132}

\bibitem{shenapde2017}
\bysame, \emph{Boundary estimates in elliptic homogenization}, Anal. PDE
  \textbf{10} (2017), no.~3, 653--694. \MR{3641883}

\bibitem{ShenBook18}
\bysame, \emph{Periodic homogenization of elliptic systems}, Operator Theory:
  Advances and Applications, vol. 269, Birkh\"{a}user/Springer, Cham, 2018,
  Advances in Partial Differential Equations (Basel). \MR{3838419}

\bibitem{shen-zhuge2016}
Z.~Shen and J.~Zhuge, \emph{Approximate correctors and convergence rates in
  almost-periodic homogenization}, J. Math. Pures Appl. (9) \textbf{110}
  (2018), 187--238. \MR{3744924}

\bibitem{tartar1985}
L.~Tartar, \emph{Estimations fines des coefficients homog\'{e}n\'{e}is\'{e}s},
  Ennio {D}e {G}iorgi colloquium ({P}aris, 1983), Res. Notes in Math., vol.
  125, Pitman, Boston, MA, 1985, pp.~168--187. \MR{909716}

\bibitem{Tbook09}
\bysame, \emph{The general theory of homogenization}, Lecture Notes of the
  Unione Matematica Italiana, vol.~7, Springer-Verlag, Berlin; UMI, Bologna,
  2009, A personalized introduction. \MR{2582099}

\bibitem{guenean2018}
N.~Wellander, S.~Guenneau, and E.~Cherkaev, \emph{Two-scale cut-and-projection
  convergence; homogenization of quasiperiodic structures}, Math. Methods Appl.
  Sci. \textbf{41} (2018), no.~3, 1101--1106. \MR{3762336}

\bibitem{xuniucpde2020}
Y.~Xu and W.~Niu, \emph{Homogenization of elliptic systems with stratified
  structure revisited}, Comm. Partial Differential Equations \textbf{45}
  (2020), no.~7, 655--689. \MR{4120920}

\end{thebibliography}

\medskip
\end{document}